\documentclass[final]{siamart190516} %[final,leqno,letterpaper]

\title{Mixed Precision Iterative Refinement with Sparse Approximate Inverse Preconditioning%\thanks{Version of \cred{\today}}
}
\author{Erin Carson\thanks{
Faculty of Mathematics and Physics, Charles University. Both authors were supported by Charles University PRIMUS project no. PRIMUS/19/SCI/11. The first author was additionally supported by Charles University
Research program no. UNCE/SCI/023 and by the Exascale Computing Project (17-SC-20-SC), a collaborative effort of the U.S. Department of Energy Office of Science and the National Nuclear Security Administration.
} \and Noaman Khan\footnotemark[2]}

\usepackage{mathtools}
\usepackage{amssymb}
\usepackage{graphicx}
\usepackage{epstopdf}
\usepackage{color}
\usepackage{graphicx}
\usepackage{lipsum}
\usepackage{graphicx}
\usepackage{subcaption}
\usepackage{hyperref}
\usepackage{algpseudocode}
\usepackage{algorithmicx}
\usepackage[T1]{fontenc}
\usepackage{hyphenat}
\usepackage[section]{placeins}
\usepackage{multirow,array,graphics}

\newcommand{\pl}{\hspace{-2pt}+\hspace{-2pt}}
\newcommand{\Mhat}{\widehat{M}}
\newcommand{\cond}{\text{cond}}
\newcommand{\mhat}{\widehat{m}}
\newcommand{\xhat}{\widehat{x}}
\newcommand{\dhat}{\widehat{d}}
\newcommand{\vhat}{\widehat{v}}
\newcommand{\yhat}{\widehat{y}}
\newcommand{\what}{\widehat{w}}
\newcommand{\zhat}{\widehat{z}}
\newcommand{\rhat}{\widehat{r}}

\newcommand{\tg}{\tilde{\gamma}}
\newcommand{\Pc}{P}
\newcommand{\Pch}{\widehat{\Pc}}

\newcommand{\cred}[1]{{\color{red}{#1}}}

\begin{document}

\maketitle

\renewcommand{\thefootnote}{\fnsymbol{footnote}}

\begin{abstract}
With the commercial availability of mixed precision hardware, mixed precision GMRES-based iterative refinement schemes have emerged as popular approaches for solving sparse linear systems. Existing analyses of these approaches, however, are based on using full LU factorizations to construct preconditioners for use within GMRES in each refinement step. In practical applications, inexact preconditioning techniques, such as incomplete LU or sparse approximate inverses, are often used for performance reasons. 

In this work, we investigate the use of sparse approximate inverse preconditioners based on Frobenius norm minimization within GMRES-based iterative refinement. We analyze the computation of sparse approximate inverses in finite precision and derive constraints under which user-specified stopping criteria will be satisfied. We then analyze the behavior of and convergence constraints for a five-precision GMRES-based iterative refinement scheme that uses sparse approximate inverse preconditioning, which we call SPAI-GMRES-IR. Our numerical experiments confirm the theoretical analysis and illustrate the resulting  tradeoffs between preconditioner sparsity and GMRES-IR convergence rate.

\end{abstract}



\section{Introduction}
We consider the problem of solving linear systems $Ax=b$ with a nonsingular $n\times n$ matrix $A$.  With the recent emergence of commercially available mixed precision hardware, there has been a renewed interest in the development of mixed precision algorithms for numerical linear algebra. The benefit of using low precision is that it is much faster and more energy efficient. Compared to double precision, for example, using half precision tensor cores available on modern NVIDIA GPUs leads to a theoretical $16\times$ improvement in performance. Low precision can also reduce communication and memory bottlenecks, since we need to move and store fewer bits. The downside is that with fewer bits, we have less accuracy and a smaller range of representable numbers. The size in bits, range, and unit roundoff $u$ are given for four IEEE precisions in Table \ref{tab:ieeeprecisions}.

\begin{table}[H]
\centering
\caption{Parameters for IEEE floating point precisions. The range denotes the order of magnitude of the largest and smallest positive normalized floating point numbers. }
\label{tab:ieeeprecisions}
\begin{tabular}{l|l|l|l}
Type & Size & Range & Unit Roundoff $u$ \\ \hline
half & 16 bits & $10^{\pm 5}$ & $2^{-11} \approx 4.9\cdot 10^{-4}$ \\
single & 32 bits & $10^{\pm 38}$ & $2^{-24} \approx 6.0\cdot 10^{-8}$ \\
double & 64 bits & $10^{\pm 308}$ & $2^{-53} \approx 1.1\cdot 10^{-16}$ \\
quad & 128 bits & $10^{\pm 4932}$ & $2^{-113} \approx 9.6\cdot 10^{-35}$
\end{tabular}
\end{table}

The goal is thus to design mixed precision approaches which use lower precision in select computations and higher precision in others, in such a way that both (1) performance is improved and (2) the attainable accuracy remains sufficient; see \cite{abdelfattah2021survey and \cite{higham2022mixed} for recent surveys}. For the problem of solving linear systems, mixed precision iterative refinement has been the focus of renewed attention. The general idea behind iterative refinement is that after computing an initial solution $x_0$ to $Ax=b$, we can refine this solution by repeatedly computing the residual $r_i = b-Ax_i$, solving $Ad_i = r_i$, and then updating the approximate solution $x_{i+1} = x_i + d_i$ until some convergence criterion is satisfied. Typically, an LU factorization of $A$ is used to solve for $x_0$, and then is  reused in each iteration to solve for the correction $d_i$. We call this LU-based variant ``standard IR'' or SIR for short. 

There is a long history of using mixed precision within iterative refinement. What we call ``traditional'' iterative refinement involves computing the residuals in double the working precision $u$, which was used already by Wilkinson in 1948 \cite{wilk48a}, and was analyzed by Wilkinson \cite{wilk63} and Moler \cite{mole67}. Fixed precision iterative refinement, in which all computations are performed in precision $u$, was analyzed by Jankowski and Wo\'{z}niakowski \cite{jawo77} and Skeel \cite{skee80}. Motivated by the faster speed of single precision versus double, in the early 2000s, there was a flurry of work in using lower precision in the computation of the LU factorization, which is the most expensive part of the computation, and the working precision elsewhere; see, for example, \cite{lllk06}, \cite[Section 9]{aadg16}. See \cite[Table 1.1]{carson2018accelerating} for a summary. 

The inclusion of half precision in modern GPUs inspired the development of iterative refinement schemes that use three or more hardware precisions. In \cite{carson2018accelerating}, the authors define an iterative refinement scheme which uses three potentially different precisions: $u_f$ for the factorization, $u$ for the working precision, and $u_r$ for the residual computation. To allow for general solvers for the correction term $d_i$, the authors also introduce a fourth precision $u_s$, which denotes the ``effective precision'' of the solve. For standard iterative refinement, the effective precision of the solve (which involves triangular solves with the LU factors computed in precision $u_f$), we have $u_s = u_f$. Then assuming that $u_f \geq u$ and $u_r \leq u^2$, the relative forward and backward errors in standard iterative refinement will converge to level $u$ when $\kappa_\infty(A)\leq u_f^{-1}$. 

The reason for introducing this effective precision becomes clear when we use a more accurate solver. The GMRES-based iterative refinement scheme (GMRES-IR) introduced in \cite{carson2017new} uses GMRES left preconditioned by the computed LU factors in order to solve for $d_i$ in each refinement step. Assuming that GMRES is executed in working precision $u$, with the matrix-vector products with the preconditioned matrix computed in precision $u^2$, GMRES-IR is guaranteed to give forward and backward errors to the working precision for more ill-conditioned systems than standard iterative refinement. Again under the assumption that $u_f \geq u$ and $u_r \leq u^2$, the relative forward and backward errors in GMRES-IR will converge to level $u$ when $\kappa_\infty(A)\leq u^{-1/2}u_f^{-1}$.

The requirement that the preconditioned matrix is applied in double the working precision within GMRES is unattractive from a performance perspective. In \cite{h:21}, the authors introduce and analyze a five-precision variant of GMRES-IR. In addition to the working precision $u$, factorization precision $u_f$, and residual precision $u_r$ for the refinement scheme, they also add precisions $u_g$ for the working precision within GMRES and $u_p$ for the precision in which the preconditioned matrix is applied to a vector within GMRES. Setting $u=u_g=u_p$ is a variant commonly used in practice. The cost is that this variant is only guaranteed to converge for a smaller range of condition numbers; again assuming that $u_f \geq u$ and $u_r \leq u^2$, one obtains relative forward and backward errors to the level of the working precision for matrices with $\kappa_\infty(A) \leq u^{-1/3}u_f^{-2/3}$, although this restriction is likely overly pessimistic in practice.

Existing analyses of GMRES-IR makes the assumption that an LU factorization is computed for use as a left preconditioner in GMRES in each refinement step. This may be undesirable from a performance perspective in the case of sparse $A$. Even if $A$ is very sparse, its LU factors may have considerable fill-in. In practice, inexact preconditioners, such as incomplete LU factorization or sparse approximate inverses are often used. 
SPAI preconditioners have the advantage that their construction is highly parallelizable; in theory, each column of the approximate inverse can be computed independently. Further, their application to a vector involves only a sparse matrix-vector multiplication, rather than the two triangular solves required for applying LU factors.

In this work, we explore the potential for the use of sparse approximate inverse preconditioners within five-precision GMRES-IR. In Section \ref{sec:background} we give a brief background on sparse approximate inverse preconditioning and mixed precision iterative refinement. In Section \ref{sec:analysis}, we analyze the error in constructing the sparse approximate inverse in finite precision and then analyze the behavior of five-precision GMRES-based iterative refinement schemes based on sparse approximate inverses, giving theoretical results on the attainable accuracy and criteria for convergence of the refinement scheme. In Section \ref{sec:experiments} we present numerical experiments and in Section \ref{sec:conclusion} we discuss open problems and future work.

\section{Background and related work}
\label{sec:background}
\subsection{Notation}
We first introduce notation which will be used in the remainder of the text. Of particular importance will be various condition numbers. For a given norm $p$, a matrix $A$, and a vector $x$, we define
\begin{equation*}
    \kappa_p(A) = \Vert A^{-1}\Vert_p \Vert A\Vert_p,\quad
    \cond_p(A) = \Vert |A^{-1}||A|\Vert_p,\quad 
    \cond_p(A,x) = \frac{\Vert |A^{-1}||A||x|\Vert_p}{\Vert x \Vert_p},
\end{equation*}
where $|A|=(|a_{ij}|)$. In the case that $p$ is not specified, the infinity norm should be assumed. We will use $u$'s to denote unit roundoffs, where subscripts on $u$ will distinguish between various precisions. For rounding error analysis, we will frequently use the notation 
\begin{equation*}
    \gamma_k = \frac{ku}{1-ku}, \quad \tilde\gamma_k=\frac{cku}{1-cku},
\end{equation*}
where $c$ is a small constant independent of problem dimension.  A superscript on $\gamma$ indicates that $u$ has that superscript as a subscript, e.g., $\gamma_k^f = ku_f/(1-ku_f)$. Quantities computed in finite precision will be denoted with hats. In order to make clear the distinction between iterative refinement (the outer solver) and GMRES (the inner solver), we will always use the word ``steps'' when referring to iterative refinement and ``iterations'' when referring to GMRES.

\subsection{Sparse approximate inverse preconditioners}

The idea behind sparse approximate inverse (SPAI) preconditioning is to explicitly construct a matrix $M\approx A^{-1}$. This has advantages within Krylov subspace methods since the application of the preconditioner involves only a matrix-vector product, rather than, e.g., the triangular solves involved when using LU-based preconditioners. There are many possible algorithms for computing $M$; see \cite{benzi1999comparative, benzi2002preconditioning} for a survey and references. 

For computing a sparse approximate inverse in unfactored form (i.e., a single matrix $M$), a popular approach is based on Frobenius norm minimization, in which $M$ is computed as the solution to $\min_{\mathcal{J}\in\mathcal{S}} \Vert I-AM\Vert_F$, where $\mathcal{J}\in \mathbb{B}^{n\times n}$ is a prescribed binary sparsity pattern in the set of all possible binary sparsity patterns $\mathcal{S}\in \mathbb{B}^{n\times n}$. The key feature of this minimization problem is that it can be decoupled as
\begin{equation}
\min_{\mathcal{J}\in\mathcal{S}} \Vert I-AM\Vert_F^2 = \sum_{k=1}^n \min_{\mathcal{J}_k\in\mathcal{S}_k} \Vert e_k-Am_k\Vert_2^2,
\label{eq:fronormmin}
\end{equation}
where $\mathcal{J}_k$, $m_k$, and $e_k$ are the $k$th columns of $\mathcal{J}$, $M$, and $I$, respectively. 
The computation of $M$ thus reduces to solving a linear least squares problem for each column $m_k$ of $M$. A primary advantage in terms of performance is that these linear least squares problems are independent, and thus can in principle be solved in parallel.

Early work based on this technique used a fixed sparsity pattern $\mathcal{J}$. The set $\mathcal{J}_k$ is called the extraction set, giving the column indices of $A$ that are relevant for solving for a column $m_k$. The nonzero rows of the submatrix $A(:, \mathcal{J}_k)$ are denoted by the so-called ``shadow'' of $\mathcal{J}_k$, 
\[
\mathcal{I}_k = \left\{ i\in\{1,\ldots, n\}: \sum_{j\in \mathcal{J}_k } |a_{ij}|\neq 0\right\},
\]
where $a_{ij}$ is the $(i,j)$ entry of $A$. Thus each term in the summation on the right in \eqref{eq:fronormmin} can be reduced to 
\begin{equation}
\min_{\mathcal{J}(\bar{m}_k) = \mathcal{J}_k} \Vert \bar{e}_k - \bar{A}_k \bar{m}_k \Vert_2,
\label{redprob}
\end{equation}
where $\bar{A}_k = A(\mathcal{I}_k, \mathcal{J}_k)\in \mathbb{R}^{|\mathcal{I}_k|,|\mathcal{J}_k|}$, $\bar{m}_k = m_k(\mathcal{J}_k)\in \mathbb{R}^{|\mathcal{J}_k|}$, $\bar{e}_k = e_k(\mathcal{I}_k)\in\mathbb{R}^{|\mathcal{I}_k|}$, and $\mathcal{J}(\bar{m}_k)$ is the binary sparsity pattern of $\bar{m}_k$. Because the resulting least squares problems are small, they can be solved directly, for example, via QR factorization.

A deficiency of this approach is that it is difficult in general to select a priori a sparsity pattern which will produce an effective preconditioner. Common choices include the sparsity pattern of $A$, $A^T$, or a power of a sparsified $A$, although these are not in general guaranteed to produce effective preconditioners. To overcome this, many authors have developed iterative approaches, in which one starts with an initial sparsity structure and adds nonzeros to this pattern until the constraint $\Vert e_k - Am_k\Vert_2\leq \varepsilon$ is satisfied for some threshold $\varepsilon$ or a maximum number of nonzeros has been reached. 
%,  
For algorithms of this type, see, e.g., the work of Cosgrove et al. \cite{cosgrove1992approximate}, Grote and Huckle \cite{grote1997parallel}, and Gould and Scott \cite{gould1998sparse}.

The algorithm of Grote and Huckle \cite{grote1997parallel} is among the most successful and commonly used approaches for computing an SPAI preconditioner \cite{benzi1999comparative}, and is the focus of our present work. Motivated by the difficulty of selecting a sparsity pattern that results in a good preconditioner a priori, the idea in \cite{grote1997parallel} was to develop an adaptive technique which captures the most important nonzero indices dynamically. 

A particular variant of Grote and Huckle's approach, which appears in \cite[Algorithm 4]{sedlacek2012sparse} is shown in Algorithm \ref{alg:spai}. 
The algorithm takes as input the matrix $A$, 
an initial binary  sparsity pattern $\mathcal{J}$, the convergence tolerance $\varepsilon$, $\alpha$, which gives the maximum number of iterations for each column, and $\beta$, which gives the maximum number of nonzeros that are added to the pattern in each iteration.

Given an initial sparsity pattern $\mathcal{J}$, for each column, the algorithm solves the least squares problem \eqref{redprob} and computes the residual $\bar{s}_k$ (lines \ref{solv1}-\ref{solvn}). If the 2-norm of the residual is below the threshold $\varepsilon$, we are finished for this column. Otherwise, we proceed with adding entries to $\mathcal{J}$. 

The nonzeros entries in $\bar{s}_k$ are given by the index set $\mathcal{L}_k$ in line \ref{ell}. For every entry $\ell$ of $\mathcal{L}_k$, we can define the set $\mathcal{N}_\ell$ which contains the column indices of the nonzero entries in the $\ell$th row of $A$ which are not already in $\mathcal{J}_k$. The union of the sets $\mathcal{N}_\ell$ gives the set $\tilde{\mathcal{J}}_k$, which is the set of potential indices to be added to $\mathcal{J}_k$. We want to select only a subset of the ``most important'' of these indicies to add. 

There are various ways of determining which indices are most important. The approach of Grote and Huckle is to consider a univariate minimization problem, through which the quantity $\rho_{jk}$ computed in line \ref{rhojk} gives a measure of the 2-norm of the new residual if index $j$ is added to $\mathcal{J}_k$. A well-known heuristic (see, e.g., \cite{barnard1999mpi}) is to mark indices as ``acceptable'' if their $\rho_{jk}$ is less than the arithmetic mean $\bar{\rho}_k$ over all $j$. In each of the $\alpha$ iterations, we add up to $\beta$ of the best (smallest $\rho_{jk}$) acceptable indices (lines \ref{add1}-\ref{addn}).

We note that the QR factorization in line \ref{qrfact} does not need to be recomputed entirely in each step, but can instead be updated using the QR factorization computed in the previous step and the entries added to $\bar{A}_k$; see \cite[Eqns.~(14) and (15)]{grote1997parallel}. A benefit of the adaptive approach of Grote and Huckle is that the algorithm is robust with respect to choice of initial sparsity pattern $\mathcal{J}$. A common choice, also recommended in \cite{grote1997parallel}, is to take the initial $\mathcal{J}$ to be the identity matrix. Typical choices of the other parameters are $\varepsilon\in [0.1,0.5]$, $\alpha\in\{1,\ldots,5\}$, and $\beta\in\{3,\ldots,8\}$ \cite[Section 3.1.3]{sedlacek2012sparse}.

Despite that each column can be computed in parallel, constructing an SPAI preconditioner is often costly, especially for large-scale problems; see, e.g., \cite{barnard1999mpi, benzi1999comparative, chow2001parallel, gao2021thread}. The memory requirements for SPAI scale quadratically and the computational cost cubically in the number of nonzeros per row \cite{gao2021thread}. Thus the use of low precision for storage and computation has the potential to significantly reduce this cost. 
The development of efficient sparse approximate inverse computations for modern hardware, especially GPUs, has been the subject of much recent work; see, e.g., \cite{gao2017gpu,lukash2012sparse,dehnavi2012parallel,he2020efficient}. 

We note that there are also techniques based on incomplete biconjugation which can be used to produce sparse triangular factors, so that the approximate inverse is a product of two matrices approximating the LU factors of $A$; see, e.g., \cite{benzi1998sparse}. Here we do not consider this approach nor the many other variants of sparse approximate inverses, including factorized sparse approximate inverses \cite{huckle2003factorized}, modified sparse approximate inverses \cite{huckle2007frobenius}, or incomplete sparse approximate inverses \cite{anzt2018incomplete}. We note that these could be interesting venues for extending the present work.

\begin{algorithm}[h]
\caption{Variant of sparse approximate inverse (SPAI) construction (\cite[Algorithm 4]{sedlacek2012sparse})}\label{alg:spai}
\begin{algorithmic}[1]
\Require {$A \in\mathbb{R}^{n\times n}$, $\mathcal{J} \in \mathbb{B}^{n\times n}$, $\alpha \geq 0$, $\beta \geq 0$, $\varepsilon >0$}
\Ensure {Right preconditioner $M \approx A^{-1}$, $M\in\mathbb{R}^{n\times n}$}
\For{$k=1$ to $n$}
    \State {$e_k=I(:,k)$}
    \State {$\mathcal{J}_k=\mathcal{J}(:,k)$}
    \For{step $=0$ to $\alpha$}
        \State {$\mathcal{I}_k=\bigg\{i \in \left\{1, \ldots ,n\right\}\colon \sum_{j \in \mathcal{J}_k} \left|a_{ij}\right| \neq 0 \bigg\}$} \label{solv1}
        \State {$\bar{A}_{k}=A(\mathcal{I}_k,\mathcal{J}_k)$}
        \State {$\bar{e}_{k}=e_k(\mathcal{I}_k)$}
        \State {Compute QR factorization $\bar{A}_{k} = \bar{Q} \bar{R}$} \label{qrfact}
        \State {$\bar{m}_{k}=\bar{R}^{-1}\bar{Q}^{T}\bar{e}_{k}$}
        \State {$\bar{s}_{k}= \bar{A}_k\bar{m}_{k}-\bar{e}_{k}$} \label{solvn}
        \If{$\Vert \bar{s}_{k}\Vert_{2}$ $\leq$ $\varepsilon$}
            \State {break}
        \EndIf
        \State {$\mathcal{L}_k=\mathcal{I}_k$ $\cup$ $\left\{ k \right\}$} \label{ell}
        \For {$\ell$ $\in$ $\mathcal{L}_k$}
            \State {$\mathcal{N}_\ell=\{j \colon {a}_{\ell j} \neq 0 \}$} \label{enn}
        \EndFor
        \State {$\tilde{\mathcal{J}}_k =\bigcup_{\ell \in \mathcal{L}_k} \mathcal{N}_\ell$}
        \State {$\tilde{\rho}_{k}=0$}
        \For {$j$ $\in$ $\tilde{\mathcal{J}}_k$}
            \State {${\rho}_{jk}=\left(\Vert \bar{s}_{k}\Vert^{2}_{2}-\frac{[\bar{r}^{T}_{k}A_j(\mathcal{I}_k)]^2}{\Vert A_j(\mathcal{I}_k)\Vert^{2}_{2}}\right)^\frac{1}{2}$} \label{rhojk}
            \State {$\tilde{\rho}_{k}=\tilde{\rho}_{k}$ + $\tilde{\rho}_{jk}$}
        \EndFor
        \State {$\tilde{\rho}_{k}=\frac{\tilde{\rho}_{k}}{|\tilde{\mathcal{J}}_{k}|}$}
        \For{idx $=1$ to $\beta$}
            \State {$j={argmin}_{j \in \tilde{\mathcal{J}}_{k}}$ $\rho_{jk}$} \label{add1}
            \State {$\mathcal{J}_k={\mathcal{J}}_{k}$ $\cup$ $\left\{j \colon \rho_{jk} \leq \tilde{\rho}_{k}\right\}$}
            \State{$\tilde{\mathcal{J}}_{k}=\tilde{\mathcal{J}}_{k}$\textbackslash $\left\{j\right\}$} \label{addn}
        \EndFor
    \EndFor
    \State {$m_k(\mathcal{J}_k)=\bar{m}_{k}$}
\EndFor
\end{algorithmic}
\end{algorithm}

\subsection{Mixed precision iterative refinement}

In Algorithm \ref{alg:genir} we present a general three-precision iterative refinement scheme. There are three explicit hardware precisions involved: $u_f$, in which the initial approximate solution $x_0$ is computed, $u_r$, in which the residual is computed, and $u$, the working precision in which all other computations are performed and all quantities are stored.

\begin{algorithm}[h]
\caption{General Iterative Refinement Scheme}\label{alg:genir}
\begin{algorithmic}[1]
\Require {$n$-by-$n$ nonsingular matrix $A$ and length-$n$ right-hand side vector $b$.}
\Ensure {Approximate solution $\xhat$ to $Ax=b$.}
\State {Compute LU factorization of $A$ in precision $u_f$.}
\State {Solve $Ax_0=b$ in precision $u_f$; store $x_0$ in precision $u$.}
\For {$i = 0,1\ldots$ until convergence}
\State {Compute $r_i=b-Ax_i$ in precision $u_r$; store $r_i$ in precision $u$.}
\State {Solve $Ad_i = r_i$ at effective precision $u_s$; store $d_i$ in precision $u$.} \label{solve}
\State {Update $x_{i+1}=x_i+d_i$ in precision $u$.}
\EndFor
\end{algorithmic}
\end{algorithm}

The effective precision $u_s$ for the solve in line \ref{solve} of Algorithm \ref{alg:genir} depends on the precisions and the solver used. In particular, $u_s$ is defined such that the solver satisfies the conditions:
\begin{align}
    \dhat_i &=(1+u_s E_i) d_i, \quad u_s \Vert E_i\Vert_\infty <1, \quad\text{and}\label{ass1}\\
    \Vert \rhat_i - A\dhat_i \Vert_\infty &\leq u_s (c_1 \Vert A\Vert_\infty \Vert \dhat_i \Vert_\infty + c_2 \Vert \rhat_i \Vert_\infty), \label{ass2}
\end{align}
where $E_i$, $c_1$, and $c_2$ are functions of $n$, $A$, $\rhat_i$, and $u_s$. These two conditions are used in analyzing the relative forward error and normwise backward errors, respectively. For simplicity, we will only consider normwise forward and backward errors in this work, but direct the reader to \cite[Section 5]{carson2018accelerating} for detailed treatment of componentwise backward error. 
Given a choice of solver with a particular $u_s$ along with precisions $u_f$, $u$, and $u_r$, the authors in \cite{carson2018accelerating} prove the following results about the behavior of the resulting iterative refinement scheme.

\begin{theorem} \cite[Corollary 3.3]{carson2018accelerating}\label{thm:fe}
Let Algorithm \ref{alg:genir} be applied to a linear system $Ax=b$ where $A$ is $n\times n$ and nonsingular with at most $q$ nonzeros per row and assume the solver used satisfies \eqref{ass1}. Then for refinement step $i\geq 0$, as long as 
\begin{equation}
\phi_i = 2u_s \min (\normalfont{\cond(A)}, \kappa_\infty(A)\mu_i)+u_s \Vert E_i \Vert_\infty 
    \label{phi}
\end{equation}
is less than 1, the forward error is reduced on the ith step by a factor of approximately $\phi_i$ until an iterate $\xhat$ is obtained for which
\begin{equation*}
    \frac{\Vert x-\xhat \Vert_\infty}{\Vert x\Vert_\infty} \lesssim 4qu_r \normalfont{\cond(A,x)} + u. 
\end{equation*}
\end{theorem}

The quantity $\mu_i$ in the above Theorem is defined such that 
\[
\Vert A(x-\xhat_i)\Vert_\infty = \mu_i \Vert A\Vert_\infty \Vert x-\xhat_i\Vert_\infty,
\]
where $\mu_i$ satisfies $\kappa_\infty(A)^{-1} \leq \mu_i \leq 1$. The insight in \cite{carson2017new} is that $\mu_i$ will be very small at the beginning of the refinement process, and will only grow close to its maximum once the process is close to convergence. Thus the size of the quantity $u_s\Vert E_i\Vert_\infty$ will usually dictate convergence in practice. 

\begin{theorem} \cite[Corollary 4.2]{carson2018accelerating}\label{thm:be}
Let Algorithm \ref{alg:genir} be applied to a linear system $Ax=b$ where $A$ is $n\times n$ and nonsingular with at most $q$ nonzeros per row and assume the solver used satisfies \eqref{ass1}. Then for refinement step $i\geq 0$, as long as 
\begin{equation}
\phi_i = (c_1 \kappa_\infty(A) + c_2)u_s 
    \label{phibe}
\end{equation}
is less than 1, the residual is reduced in each step by a factor of approximately $\phi_i$ until an iterate $\xhat$ is obtained for which
\begin{equation*}
    \Vert b - A\xhat\Vert_\infty \lesssim qu(\Vert b\Vert_\infty + \Vert A\Vert_\infty \Vert \xhat\Vert_\infty),
\end{equation*}
which indicates that $\xhat$ is a backward stable solution to the working precision. 
\end{theorem}

Thus the quantities $u_s$, $\Vert E_i\Vert_\infty$, $c_1$, and $c_2$, which come from the particular solver, completely determine the convergence conditions for the forward and backward errors in iterative refinement. 
The authors in \cite{carson2018accelerating} bound these quantities for two different solvers: 1) triangular solves using the computed LU factors (SIR), and 2) left-preconditioned GMRES, where the computed LU factors are used as preconditioners and the preconditioners and preconditioned coefficient matrix are applied to vectors in double the working precision (GMRES-IR). The authors in \cite{h:21} recently extended this analysis to cover the more general case where the solver is GMRES left-preconditioned by the computed LU factors run in a working precision $u_g$, where the preconditioners and preconditioned coefficient matrix are applied to vectors in a precision $u_p$. The variant of GMRES-IR in which $u=u_g=u_p$ is attractive from a performance standpoint and is commonly used in practice. 

We summarize the size of these quantities in Table \ref{tab:boundquantities}. Here, $f$ denotes a generic function, $n$ is the problem dimension, $k$ is the maximum number of GMRES iterations in any refinement step, $\rho_n$ is the growth factor in LU, and $\tilde{A}=\widehat{U}^{-1}\widehat{L}^{-1}A$ for the computed factors $\widehat{L}$ and $\widehat{U}$.

\begin{table}[]
\centering
\caption{Summary of sizes of quantities in the assumptions \eqref{ass1}-\eqref{ass2} for various solvers as given in \cite{h:21}. In the table, $f$ denotes a generic function, $n$ is the problem dimension, $k$ is the maximum number of GMRES iterations in any refinement step, $\rho_n$ is the growth factor in LU, and $\tilde{A}=\widehat{U}^{-1}\widehat{L}^{-1}A$ for the computed factors $\widehat{L}$ and $\widehat{U}$.}
\label{tab:boundquantities}
\begin{tabular}{c|ccc}
\multicolumn{1}{l|}{} & SIR & GMRES-IR ($u_g=u, u_p=u^2$) & GMRES-IR ($u_g=u_p=u$) \\ \hline
$u_s$ & $u_f$ & $u$ & $u_g+u_p\kappa_\infty(A)$ \\
$\Vert E_i\Vert_\infty$ & $f(n,\rho_n)$ & $f(n,k,\rho_n)\kappa_\infty(\tilde{A})$ & $f(n,k,\rho_n)\kappa_\infty(\tilde{A})$ \\
$c_1$ & $f(n,\rho_n)$ & $f(n,k,\rho_n)\Vert \tilde{A}\Vert_\infty$ & $f(n,k,\rho_n)\Vert \tilde{A}\Vert_\infty$ \\
$c_2$ & $f(n,\rho_n)$ & $f(n,k,\rho_n)\kappa_\infty(A)$ & $f(n,k,\rho_n)\kappa_\infty(A)$
\end{tabular}
\end{table}

Substituting the values for $u_s$, $\Vert E_i\Vert_\infty$, $c_1$, and $c_2$ from Table \ref{tab:boundquantities} into Theorems \ref{thm:fe} and \ref{thm:be} shows that for SIR, 
the forward and backward errors will converge as long as $\kappa_\infty(A) \ll u_f^{-1}$, with a limiting value of $qu$ for the backward error and a limiting value of $qu_r\cond(A,x)+u$ for the forward error. 

GMRES-IR with $u_g=u$ and $u_p=u^2$ will have the same limiting values, but looser constraints on condition number. Using the bound $\kappa_\infty(\tilde{A})\ll \kappa_\infty(A)^2u_f^2$ (see \cite[Eqn. (8.3)]{carson2018accelerating}), for the forward error, we require $\kappa_\infty(A) \ll u^{-1/2}u_f^{-1}$. For the backward error, we can use the bound on $\Vert \tilde{A}\Vert_\infty$ from \cite[Eqn.~(3.2)]{carson2017new} to give the constraint $\kappa_\infty(A)\ll u^{-1/2}u_f^{-1/2}$; see also \cite[Table 2.1]{h:21}. Note that in the case that $u_r = u^2$, the constraint for the convergence of the backward error to the working precision also becomes $\kappa_\infty(A) \leq u^{-1/2}u_f^{-1}$ since the backward error is bounded by the forward error.

From the analysis in \cite{h:21}, for GMRES-IR with $u=u_g=u_p$, using the same bound on $\kappa_\infty(\tilde{A})$, the constraint for convergence of the forward error to the limiting value is $\kappa_\infty(A) \leq u^{-1/3}u_f^{-2/3}$ and that for the backward error is $\kappa_\infty(A) \leq u^{-1/3}u_f^{-1/3}$. As above, the constraint for convergence of the backward error is stricter, and thus we can take the constraint for the convergence of both forward and backward errors to the limiting value to be $\kappa_\infty(A) \leq u^{-1/3}u_f^{-2/3}$ in the case $u_r=u^2$. We stress that these constraints are overly strict in practice, and convergence is often observed even in cases where these conditions are violated; see, e.g., \cite{oktay2022multistage}.

In order to simplify the presentation and allow easy extension of results, we will restrict ourselves to left preconditioning, since this has been used in previous works on GMRES-based iterative refinement; see, e.g., \cite{h:21, carson2018accelerating}.
In \cite[Remark 3.5]{arioli2009using}, Arioli and Duff comment that right-preconditioned GMRES is only backward stable if the preconditioner satisfies certain constraints. For the case of sparse approximate inverses, this condition is $\Vert \Mhat^{-1} \Vert_2 / \Vert A \Vert_2 \approx 1$, where $\Mhat$ is the computed sparse approximate inverse. If we could guarantee this, then it could be shown that right-preconditioned GMRES is backward stable and we conjecture that similar bounds could be proven on the quantities in Table \ref{tab:boundquantities}. We note that we could also use right-preconditioned FGMRES with an SPAI preconditioner within GMRES-based iterative refinement, although such a result has not yet been formally proven. We leave these tasks as future work, but note that it is likely that right-preconditioning will work in practice. In the context of SPAI, right-preconditioning may be preferred when $A^{-1}$ has some dense rows, since this will yield a sparser approximate inverse \cite[pp.~843]{grote1997parallel}.

\subsection{Related Work}

Along with the resurgence of interest in mixed precision iterative refinement have come a number of works that build on and expand the work on GMRES-based iterative refinement, e.g., \cite{carson2017new, carson2018accelerating, aadg16, h:21}. In \cite{higham2019new}, Higham and Mary develop a new general preconditioning technique based on a low-rank approximation of the error $E=\widehat{U}^{-1}\widehat{L}^{-1}A-I$, where $\widehat{L}$ and $\widehat{U}$ are approximate LU factors of $A$. A randomized SVD is used to produce this low-rank approximation, and the preconditioner is applied in the context of GMRES-IR. This technique can apply to the case where the approximation in LU is due to low precision computation, as well as other sources of error including using an incomplete LU factorization and a block low-rank LU factorization. See \cite{amestoy2022,hima21} for other theoretical and experimental results on the use of block low-rank LU factorization preconditioners within iterative refinement. 

In \cite{lld2021}, Lindquist, Luszczek, and Dongarra have recently evaluated mixed precision restarted GMRES (equivalent in some sense to a variant of GMRES-IR) on a GPU-accelerated node with various preconditioners. They found that the use of incomplete LU with zero fill-in (ILU(0)) resulted in unfavorable performance (a slowdown on average). The authors attribute this partially to the fact that sparse triangular solves have limited parallelism for the GPU to exploit. 

In the recent work \cite{amestoy2022}, Amestoy et al. extend the work on GMRES-based iterative refinement to LU factorization methods more appropriate for sparse matrices, including those based on block low-rank approximations and on static pivoting strategies. 

The present work differs from these recent approaches in that we depart from the restriction to preconditioning based on LU factorization and instead look at sparse approximate inverses as preconditioners within mixed precision GMRES-based iterative refinement. In the following section, we will first analyze the numerical properties of SPAI preconditioners computed in a given precision with given parameters, and then evaluate the constraints for convergence of GMRES-IR with SPAI preconditioners.

\section{Analysis of SPAI-based iterative refinement}
\label{sec:analysis}

Our ultimate goal is to prove conditions under which the SPAI-GMRES-IR algorithm (Algorithm \ref{alg:spaigmresir}) will converge and to determine how the precision $u_f$ and approximation parameter $\varepsilon$ should be related. We first give bounds on the quality of the SPAI preconditioner when it is computed in finite precision, analogous to those given for the case of exact arithmetic in \cite{grote1997parallel}.

Suppose we want to construct a sparse approximate inverse of a nonsingular matrix $A$ for use as a left preconditioner, which we will denote $\Pc$. That is, instead of solving $Ax=b$ we will solve $\Pc Ax = \Pc b$. We will construct $\Pc$ using the algorithm of Grote and Huckle \cite{grote1997parallel} (see Algorithm \ref{alg:spai}), executed in some precision $u_f$. Algorithm \ref{alg:spai} as stated produces a right preconditioner $M$. We can obtain a left preconditioner by executing the algorithm on $A^T$ and then transposing the obtained approximate inverse, i.e., $\Pc = M^T$.
We will denote the version of $M$ computed in finite precision as $\Mhat$. 

 In contrast to how the method is often used in practice, we will not specify a maximum number of nonzeros here (i.e., in the extreme case we could allow a fully dense inverse). More important for our purposes will be the assumption that for each column $\mhat_k$ of $\Mhat$, it holds that 
\begin{equation}
\Vert e_k - A^T \mhat_k \Vert_2 \leq \varepsilon.
\label{eq:trueres}
\end{equation}

The way the SPAI algorithm works, we keep iterating, computing a QR factorization of some submatrix $\bar{A}^T$ of $A^T$ (selected by specifying a particular set of row and column indices) and solving the least squares problem $\min_{\bar{m}_k} \Vert \bar{A}^T\bar{m}_k-\bar{e}_k\Vert_2$. The vector $\bar{m}_k$ then gives the nonzero entries in column $k$ of $M$ (corresponding to the same indices as the rows of $A^T$ selected to construct $\bar{A}^T$). Row and column indices are added to the set, and the iteration is continued until $\min_{\bar{m}_k} \Vert \bar{A}^T\bar{m}_k-\bar{e}_k\Vert_2 \leq \varepsilon$, which indicates that \eqref{eq:trueres} holds.

 We assume that all computations within the SPAI algorithm are performed in a precision $u_f$.  We make the simplifying assumption that the QR factorization of $\bar{A}^T$ is computed fully, but note again that in practice one only needs to compute a QR factorization of a smaller block corresponding to the newly added indices; see, e.g., \cite[Algorithm 5]{sedlacek2012sparse}. We further assume 
that the QR factorization is computed using Householder QR (or a method with similar numerical properties, like TSQR \cite{demmel2012communication}).

 Now, assuming we impose no maximum sparsity pattern on $\Mhat$, we seek to determine under what constraint on $u_f$ we can guarantee that $\Vert \rhat_k\Vert_2 \leq \varepsilon$ can be satisfied, where $\rhat = fl_{u_f}(e_k - A^T \mhat_k)$ for the computed $\mhat_k$. 
Here we will assume that no overflow or underflow occurs during the SPAI computation in precision $u_f$. We note that to mitigate the risk of this one could use sophisticated scaling techniques such as those described in \cite{higham2019squeezing} and \cite{carson2020three}; we use a similar column scaling technique in our experiments in Section \ref{sec:experiments}. To account for finite precision error, it suffices to assume that we are on the final iteration for column $k$, after which we will quit since the condition $\Vert \rhat_k\Vert_2$ is satisfied.
Using the backward error result \cite[Theorem 20.3]{high:ASNA2}, which says that the computed solution $\mhat_k$ is the exact least squares solution to 
\begin{equation*}
    \min_{\mhat_k} \Vert (e_k+\Delta e_k) - (A^T+\Delta A^T)\mhat_k \Vert_2,
\end{equation*}
where $|\Delta A^T| \leq n\tg_{n}^f ee^T |A^T|$ and $| \Delta e_k | \leq n\tg_n^f ee^T|e_k| $, we can then bound 
\begin{equation}
    \Vert e_k - A^T \mhat_k \Vert_2 \leq n^2 \tg_n^f \left\Vert |e_k| + |A^T| |\mhat_k| \right\Vert_2 + \left(1+ n^2\tg_n^f \cond_2(A)\right)\Vert r_k\Vert_2 + O(u_f^2),
    \label{bebound}
\end{equation}
where $r_k = e_k - A^T m_k$. Then using standard rounding error analysis, we have 
\begin{equation*}
    \Vert \rhat_k \Vert_2 \leq \Vert e_k - A^T \mhat_k \Vert_2 + \Vert \Delta \rhat_k\Vert_2,
\end{equation*}
where $\Vert \Delta \rhat_k \Vert_2 \leq \tg_n^f \left\Vert |e_k| + |A^T| |\mhat_k|\right\Vert_2$.
Then 
\begin{equation}
    \Vert \rhat_k \Vert_2 \leq n^2 \tg_n^f \left\Vert |e_k| + |A^T| |\mhat_k| \right\Vert_2 + \left(1+ n^2\tg_n^f \cond_2(A)\right)\Vert r_k\Vert_2 + O(u_f^2),
    \label{eq:rhatbnd}
\end{equation}
where the $\Vert \Delta \rhat_k\Vert_2$ factor has been absorbed into the $\tg_n^f$ in the first term on the right-hand side. 

Recall that our goal is to determine under what constraint on $u_f$ we can say that $\Vert \rhat_k\Vert_2 \leq \varepsilon$ is achievable, with no restriction on the number of nonzeros in the columns $\mhat_k$.
In theory, we could keep increasing the size of the subproblem in solving for the column of the inverse until we have the full matrix $A^T$, which would mean that in exact arithmetic, $m_k$ is the exact $k$th column of the inverse of $A^T$, and thus $r_k=0$. Thus the second term on the right-hand side of \eqref{eq:rhatbnd} will not limit us from satisfying the stopping criterion $\Vert \rhat_k\Vert_2 \leq \varepsilon$, and it can thus be ignored for our purposes.
We are then left with (ignoring higher order terms)
\begin{equation*}
    \Vert \rhat_k \Vert_2 \leq n^2 \tg_n^f \left\Vert |e_k| + |A^T| |\mhat_k| \right\Vert_2.
\end{equation*}
Then in order to guarantee that we will eventually iterate to a solution with $\Vert \rhat_k\Vert_2 \leq \varepsilon$, we should have
\begin{equation}
n^2 \tg_n^f \left\Vert |e_k| + |A^T| |\mhat_k| \right\Vert_2 \leq \varepsilon.
    \label{convcond}
\end{equation}
In other words, the problem must not be so ill-conditioned relative to $u_f$ such that we can not guarantee that we incur an error less than $\varepsilon$ just computing the residual. We note that the dimensional constants here can be a large overestimate, since we in practice solve the least squares problem with the smaller submatrix $\bar{A}^T$. We can turn this into the more descriptive (but less rigorous) bound 
\[
\cond_2(A^T,m_k) \lesssim \varepsilon u_f^{-1},
\]
as well as an a priori bound of the form
\begin{equation}
\cond_2(A^T) \lesssim \varepsilon u_f^{-1}.
\label{convcondloose}
\end{equation}
We note that this is likely quite pessimistic in practice. Note that $\cond_2(A^T)$ can be considerably smaller than $\kappa_2(A)$. 
 Perhaps a more useful way to think about \eqref{convcondloose} is that with a given matrix $A$ and a chosen precision $u_f$, in order to guarantee success of the SPAI algorithm, one must set $\varepsilon$ such that 
 \[
u_f\cond_2(A^T) \leq \varepsilon.
 \]

This analysis echoes what we may intuitively think: the larger we make $\varepsilon$ (meaning the more approximate the inverse), the larger the $u_f$ we can tolerate. If we require smaller $\varepsilon$, meaning we want a better quality preconditioner, then we must have a smaller $u_f$ to match.

Given that \eqref{convcond} is satisfied, we know that we can find a computed solution $\mhat_k$ such that $\Vert \rhat_k\Vert_2 \leq \varepsilon$. Writing $\rhat_k = e_k - A^T \mhat_k + \Delta \rhat_k$, we have 
\begin{align*}
\Vert e_k - A^T\mhat_k \Vert_2 &\leq \Vert \rhat_k\Vert_2 + \Vert \Delta \rhat_k\Vert_2 \\
& \leq \varepsilon + \tg_n^f \left\Vert |e_k| + |A^T| |\mhat_k| \right\Vert_2\\
&\leq 2\varepsilon,
\end{align*}
from which we can write the Frobenius norm bound
\begin{equation*}
    \Vert I -A^T\Mhat\Vert_{F} \leq 2\sqrt{n}\varepsilon.
\end{equation*}
Note that assuming exact arithmetic, e.g., using $u_f = 0$ in \eqref{bebound}, the above bound becomes $\Vert I-A^T M\Vert_F \leq \sqrt{n}\varepsilon$, which is the same bound given by Grote and Huckle \cite[Theorem 3.1]{grote1997parallel}. 
At this point we can transpose to get the desired bound for the left preconditioner $\widehat{\Pc} = \Mhat^T$, and then convert to the desired infinity norm, giving  
\begin{equation}
         \Vert I -\Pch A\Vert_{\infty} \leq 2n\varepsilon.
    \label{eq:AM1}
\end{equation}
We note that from \eqref{eq:AM1} we then also have an a priori bound on the distance from the computed $\Pch$ to the true inverse of $A$,
\begin{equation}
\Vert \Pch - A^{-1}\Vert_\infty = \Vert (\Pch A - I)A^{-1} \Vert_\infty \leq \Vert I -\Pch A\Vert_{\infty} \Vert A^{-1}\Vert_\infty \leq  2n\varepsilon \Vert A^{-1}\Vert_\infty.
    \label{eq:MtoAI}
\end{equation}

An important conclusion that we can draw from this analysis, in particular from the fact that \eqref{eq:AM1} is guaranteed as long as \eqref{convcondloose} holds, is that there is \emph{no advantage} to using a higher precision in computing the SPAI preconditioner beyond that dictated by $u_f \approx \varepsilon \cond_2^{-1}(A^T)$. Using higher precision will have very little effect, since we will end up with a preconditioner of similar quality. We will demonstrate this observation in Section \ref{sec:experiments} below. Note that this is in contrast to the case of using full LU factorization, where the quality of the factors scales with the precision used. 

\begin{figure}[H]
\centering
\includegraphics[trim={3.8cm 8cm 3.8cm 8cm},clip,width=.49\linewidth]{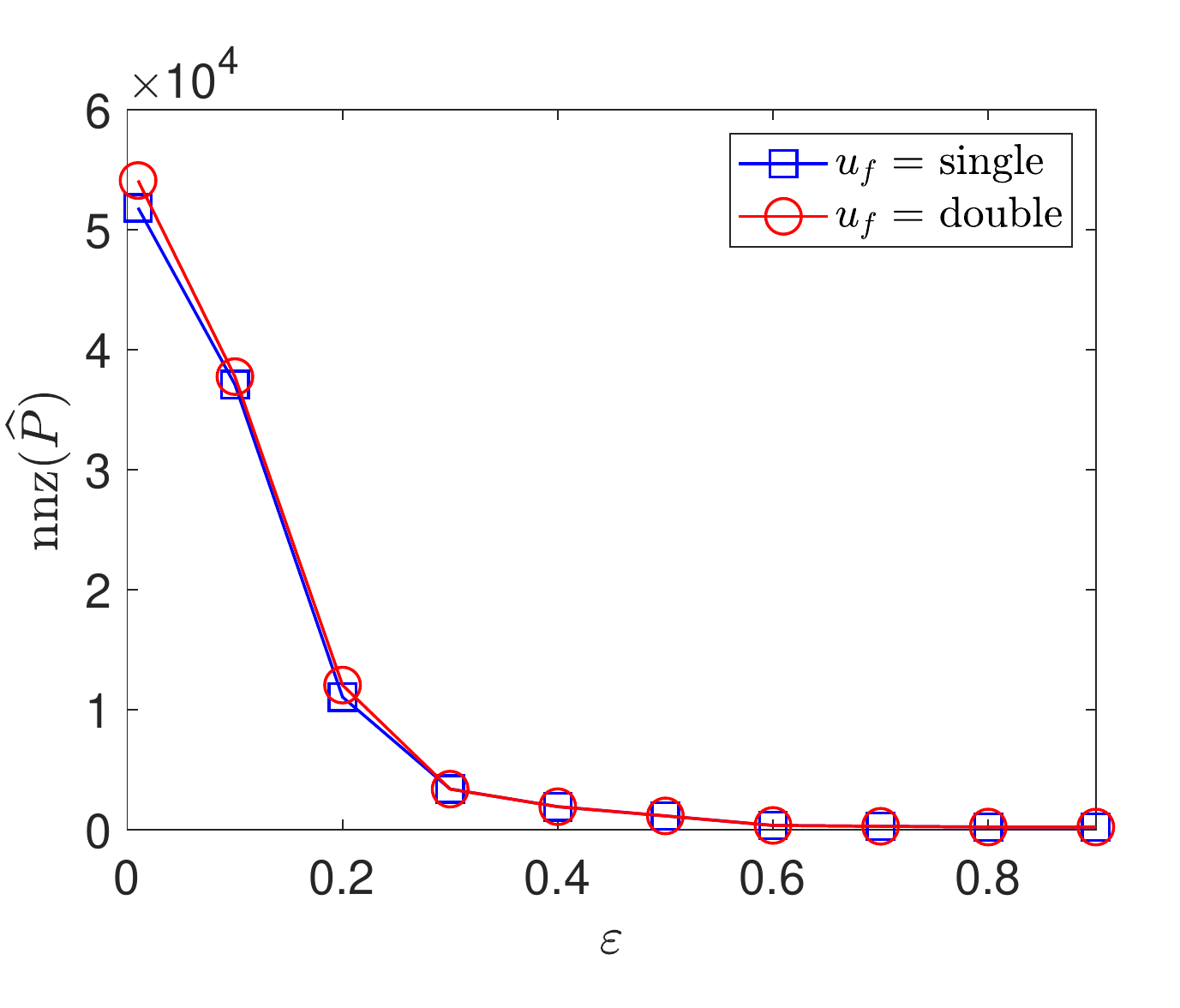} 
\includegraphics[trim={3.4cm 8cm 4cm 8cm},clip,width=.49\linewidth]{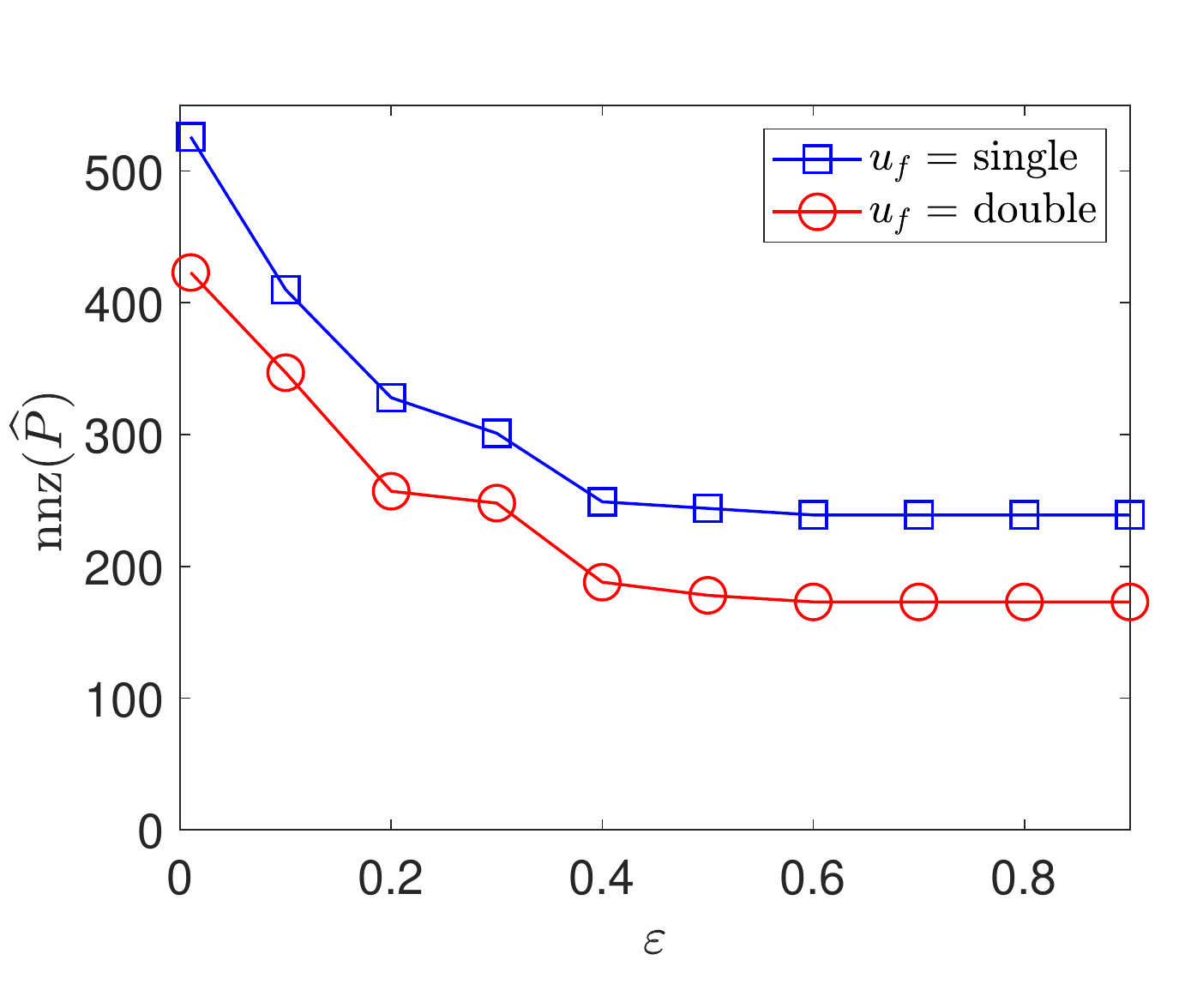} 
\caption{Number of nonzeros in the computed sparse approximate inverse preconditioner $\Pch$ when $\Pch$ is computed in single precision (blue squares) and double precision (red circles), for the matrices \texttt{saylr1} (left) and \texttt{steam3} (right), for various choices of $\varepsilon$.}
\label{fig:nnzsaylr1}
\end{figure}

Intuitively, however, we may think that using higher precision for $u_f$ might  result in a $\Pch$ with fewer (more accurate) nonzeros than a $\Pch$ computed in lower precision. This can be the case, but not necessarily. We demonstrate this in Figure \ref{fig:nnzsaylr1} for the matrices \texttt{saylr1} (left plot) and \texttt{steam3} (right plot) from SuiteSparse \cite{davis2011university}, where we plot the number of nonzeros in the resulting $\Pch$ computed using the SPAI algorithm (Algorithm \ref{alg:spai}) in single precision (blue squares) and double precision (red circles) for various $\varepsilon$ values. Within Algorithm \ref{alg:spai}, we use the identity as the initial sparsity pattern, $\beta=8$, and $\alpha = \lceil n/\beta \rceil$; this value of $\alpha$ allows the sparse approximate inverse to fill in as much as needed in order to satisfy \eqref{eq:trueres}.

For \texttt{saylr1}, there is very little relative difference between single and double precision in the resulting size of $\Pch$; if anything, single precision results in a $\Pch$ with slightly fewer nonzeros. For the matrix \texttt{steam3}, the relative difference between single and double precision is more significant; $\Pch$ computed in double precision has, on average, about 24\% fewer nonzeros than the $\Pch$ computed in single precision. Note that for both matrices, these choices of $\varepsilon$ satisfy \eqref{convcondloose} for both single and double precision; see Table \ref{tab:mats} for properties of the matrices. We will explore further examples in Section \ref{sec:ssd}.

\subsection{GMRES-based iterative refinement with SPAI}
\label{sec:spaigmresir}

Our goal is now to bound the relative error in the correction solve when the solver is GMRES left-preconditioned by the approximate inverse $\Pch$. For clarity, we present this variant, which we call SPAI-GMRES-IR, in Algorithm \ref{alg:spaigmresir}. We follow and heavily rely on the analysis of \cite{h:21}, \cite{carson2018accelerating}, and \cite{carson2017new}. Here we will analyze the general five-precision GMRES-IR variant of \cite{h:21}, in which in addition to precisions $u_f$, $u$, and $u_r$, $u_g$ denotes the working precision used within GMRES and $u_p$ represents the precision in which the preconditioner or preconditioned linear system is applied to a vector within GMRES.

\begin{algorithm}[h]
\caption{GMRES-based Iterative Refinement with SPAI Preconditioning (SPAI-GMRES-IR)}\label{alg:spaigmresir}
\begin{algorithmic}[1]
\Require {$n\times n$ nonsingular matrix $A$ and length-$n$ right-hand side vector $b$, maximum number of refinement steps $i_{max}$, GMRES convergenge tolerance $\tau$, SPAI parameter $\varepsilon$.}
\Ensure {Approximate solution $x_{i+1}$ to $Ax=b$.}
\State{Run Algorithm \ref{alg:spai} on $A^T$ with parameter $\varepsilon$ to obtain $M$; set $\Pc = M^T$.}
\State {Compute $x_0=\Pc b$ in precision $u_f$; store $x_0$ in precision $u$.}
\For {$i = 0:i_{max}-1$}
\State {Compute $r_i=b-Ax_i$ in precision $u_r$; store in precision $u$.}
\State {Solve $\Pc Ad_i = \Pc r_i$ via left-preconditioned GMRES with tolerance $\tau$ in working precision $u_g$, with matrix-vector products with $\Pc$ and $A$ computed in precision $u_p$; store $d_i$ in precision $u$.}
\State {Update $x_{i+1}=x_i+d_i$ in precision $u$.}
\EndFor
\end{algorithmic}
\end{algorithm}

We first want to bound $\kappa_\infty(\tilde{A})$, where $\tilde{A} = \Pch A$. We stress that the condition number of the preconditioned matrix will not tell us anything about the  convergence rate of GMRES (see \cite{greenbaum1996any}); here we will use it to bound the relative error after we obtain a backward error result for the preconditioned system.
We can write
\begin{align}
    \tilde{A} = \Pch A &= I-(I-\Pch A),\nonumber\\
    \tilde{A}^{-1} = A^{-1}\Pch^{-1} &= (I-(I-\Pch A))^{-1} \nonumber \\
    &\approx I+(I-\Pch A). \label{questionable}
\end{align}
Thus using \eqref{eq:AM1}, we have 
\begin{align*}
    \Vert \tilde{A} \Vert_\infty &\leq 1+2n\varepsilon,\\
    \Vert \tilde{A}^{-1} \Vert_\infty &\lesssim 
    1+2n\varepsilon,
\end{align*}
which gives
\begin{equation}
    \kappa_\infty(\tilde{A}) \lesssim (1+2n\varepsilon)^2.
    \label{kappaAtest}
\end{equation}

\begin{figure}[h]
\begin{subfigure}{0.48\textwidth}
\includegraphics[trim={3.5cm 8cm 3.5cm 8cm},clip,width=1\linewidth]{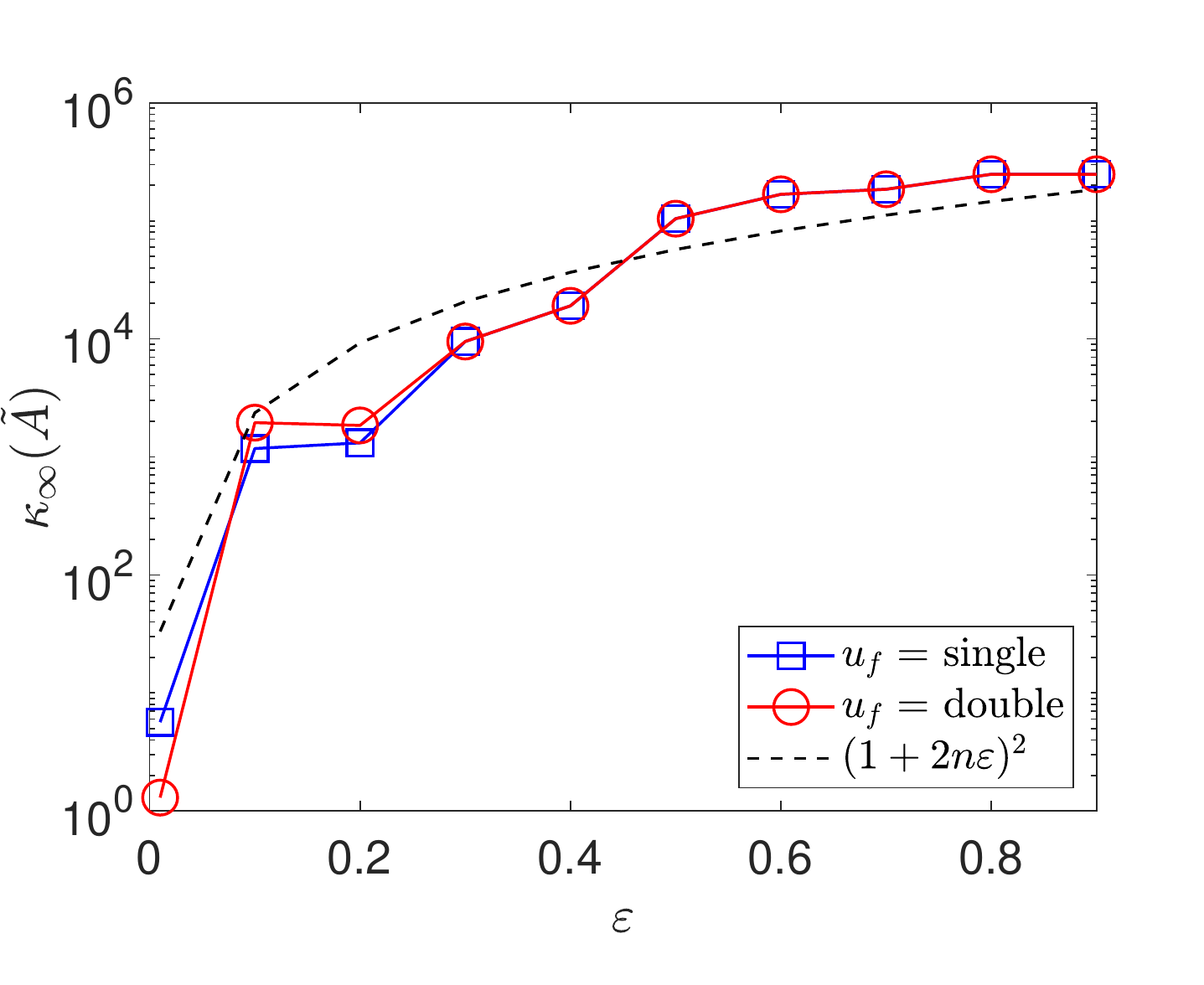} 
\caption{\texttt{saylr1}}
\label{fig:epskappasaylr1}
\end{subfigure}
\begin{subfigure}{0.48\textwidth}
\includegraphics[trim={3.5cm 8cm 3.5cm 8cm},clip,width=1\linewidth]{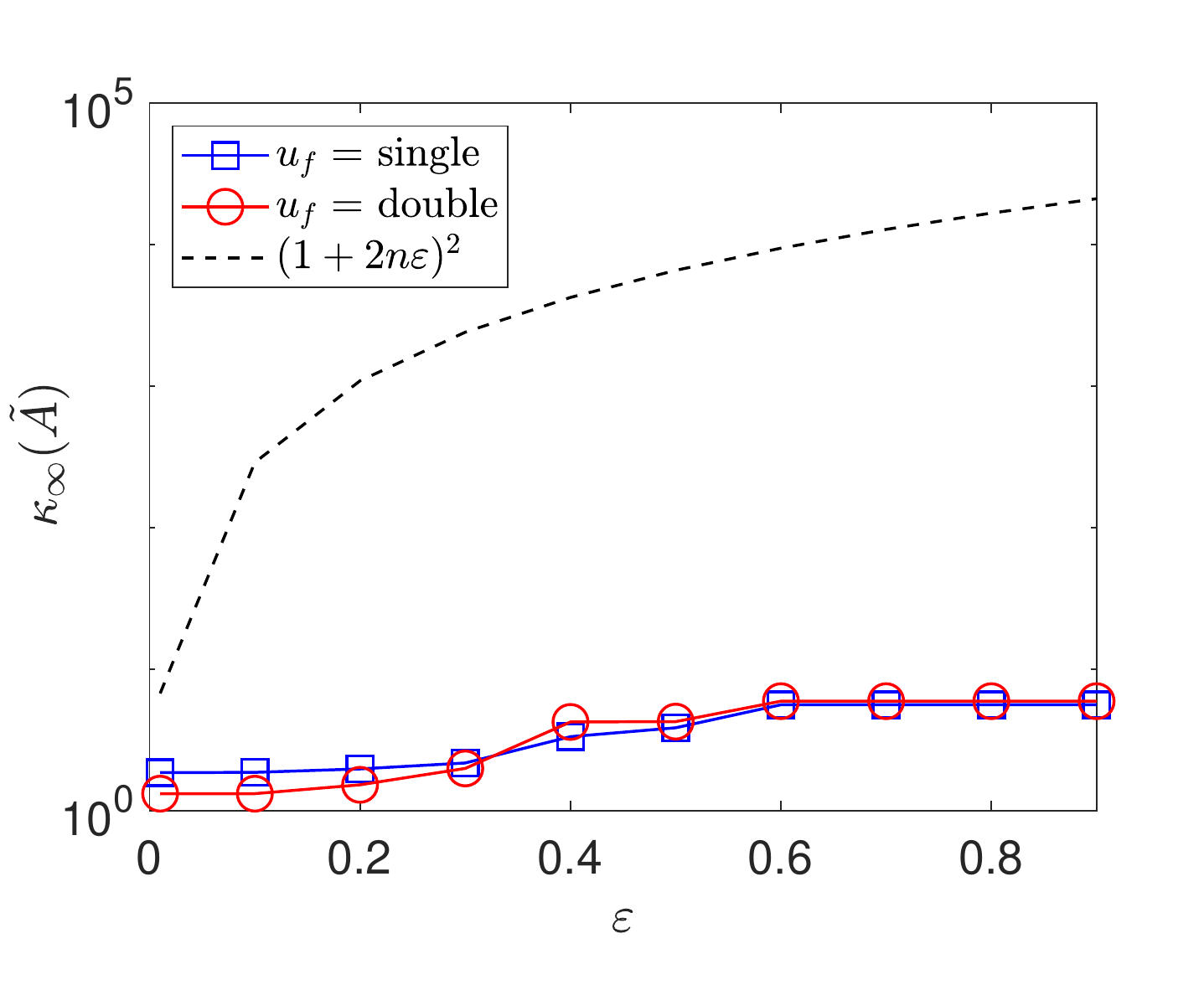}
\caption{\texttt{steam3}}
\label{fig:epskappasteam3}
\end{subfigure}
\caption{$\kappa_{\infty}(\tilde{A})$ versus $\varepsilon$ for matrices \texttt{saylr1} (left) and \texttt{steam3} (right) for SPAI preconditioners computed in single (blue squares) and double (red circles), compared with the approximation \eqref{kappaAtest}.}
\label{fig:epskappa}
\end{figure}

We note that this should not be taken as an upper bound but rather as a rough estimate, since for practical problem sizes, $2n\varepsilon$ will likely be greater than 1, and thus our dropping of higher order terms in the Taylor expansion to get \eqref{questionable} may be dubious. Although \eqref{kappaAtest} is not strictly an upper bound, we note that the growth of $\kappa_\infty(\tilde{A})$ does generally follow this bound asymptotically. In Figure \ref{fig:epskappa}, we plot the condition number $\kappa_\infty(\tilde{A})$ versus $\varepsilon$ for two matrices from SuiteSparse \cite{davis2011university}, \texttt{saylr1} (Fig. \ref{fig:epskappasaylr1}) and \texttt{steam3} (Fig. \ref{fig:epskappasteam3}), where the sparse approximate inverses are computed in single precision (blue squares) and single precision (red circles). As before, we use the identity as the initial sparsity pattern, $\beta=8$, and $\alpha = \lceil n/\beta \rceil$; see Algorithm \ref{alg:spai}.
We additionally plot $(1+2n\varepsilon)^2$, the estimate in \eqref{kappaAtest}, demonstrating that \eqref{kappaAtest} well describes the resulting condition number of the preconditioned coefficient matrix. We again note that there is little difference between single and double precision, since we expect preconditioner quality to depend on $\varepsilon$ rather than $u_f$ as long as \eqref{convcondloose} is satisfied.

Now we turn our attention to the application of the preconditioner to the right-hand side $\rhat_i$. Let the (exact) preconditioned right-hand side be $z_i = \Pch \rhat_i$. Assuming we apply $\Pch$ to $\rhat_i$ in precision $u_p$, we have 
\begin{align*}
\zhat_i &= (\Pch + \Delta \Pc)\rhat_i, \qquad |\Delta \Pc|\leq \gamma_n^p |\Pch|, \\
&=z_i + \Delta \Pc \Pch^{-1} z_i.
\end{align*}
We then have the bound 
\begin{equation}
\Vert \zhat_i - z_i \Vert_\infty \leq \gamma_n^p \Vert |\Pch| |\Pch^{-1}| \Vert_\infty \Vert z_i \Vert_\infty \leq  \gamma_n^p \kappa_\infty(\Pch)\Vert z_i \Vert_\infty
\lesssim \gamma_n^p \kappa_\infty(A)\Vert z_i \Vert_\infty,
\label{zi}
\end{equation}
under the reasonable assumption that $\kappa_\infty(\Pch) \approx \kappa_\infty(A)$.

Following \cite{carson2017new} and \cite{h:21}, we now seek to bound the error in applying $\Pch A$ to a vector $\vhat_j$ in iteration $j$ of the GMRES algorithm. We compute this via 
\begin{align*}
    (A+\delta A) \vhat_j &= \what_j, \qquad |\Delta A| \leq \gamma_n^p |A|, \\
    (\Pch + \Delta \Pch)\what_j &= \yhat_j, \qquad |\Delta \Pch|\leq \gamma_n^p |\Pch|.
\end{align*}
Then we can write the computed vector $\yhat_j$ as 
\[
\yhat_j = (\Pch + \Delta \Pch)(A+\Delta A)\vhat_j = (\Pch A + \Delta \Pch A + \Pch \Delta A)\vhat_j = (\tilde{A} + \Delta \tilde{A}')\vhat_j,
\]
where 
\begin{align*}
    \Delta \tilde{A}' & = \Delta \Pch A + \Pch \Delta A \\
    &= \Delta \Pch \Pch^{-1} \Pch A + \Pch A A^{-1} \Delta A\\
    &= \Delta \Pch \Pch^{-1} \tilde{A} + \tilde{A} A^{-1}\Delta A,
\end{align*}
for which we can write the bound
\begin{align*}
\Vert \Delta \tilde{A}' \Vert_F &\leq \gamma_n^p \big( \Vert |\Pch||\Pch^{-1}|\Vert_F + \Vert |A^{-1}||A|\Vert_F \big) \Vert \tilde{A}\Vert_F\\
&\leq \gamma_n^p \big( \kappa_F(\Pch) + \kappa_F(A)\big)\Vert \tilde{A}\Vert_F\\
&\lesssim \tilde\gamma_n^p \kappa_F(A)\Vert \tilde{A}\Vert_F.
\end{align*}
Thus
\[
\yhat_j = \tilde{A}\vhat_j + f_j,
\]
with 
\[
\Vert f_j \Vert_\infty \leq \tilde\gamma_n^p \kappa_F(A)\Vert \tilde{A}\Vert_F \Vert \vhat_j \Vert_2.
\]
Since \cite[Eqn. 3.2]{h:21} is satisfied, we can thus apply Theorem 3.1 from \cite{h:21}, which gives the result that for some iteration $k\leq n$, the left-preconditioned MGS-GMRES algorithm will produce an approximate solution $\dhat_i$ for which 
\[
(\tilde{A}+\Delta \tilde{A})\dhat_i = \zhat_i + \Delta \zhat_i, 
\]
with
\begin{equation}
 \Vert \Delta \tilde{A}\Vert_F \leq (\tg_{kn}^g + n\tg_{k^{1/2}n}^p \kappa_\infty(A)) \Vert \tilde{A}\Vert_F, \quad \Vert \Delta \zhat_i\Vert_2 \leq \tg_{kn}^g\Vert\zhat_i\Vert_2 \lesssim n^{1/2}\tg_{kn}^g \Vert z_i\Vert_\infty.
    \label{begmres}
\end{equation}
Using this result together with \eqref{zi}, we can bound the normwise relative backward error by 
\begin{equation}
\frac{\Vert z_i - \tilde{A}\dhat_i\Vert_\infty}{\Vert\tilde{A}\Vert_\infty \Vert \dhat_i\Vert_\infty + \Vert z_i \Vert_\infty} \lesssim f(n,k) \left(u_g + u_p \kappa_\infty(A)\right),
\label{nrbe}
\end{equation}
where $f(n,k)=n^3k^{1/2}$. This gives the relative forward error bound 
\begin{equation}
\frac{\Vert \dhat_i - d_i\Vert_\infty}{\Vert d_i\Vert_\infty} \lesssim f(n,k) \left(u_g + u_p \kappa_\infty(A)\right) \kappa_\infty(\tilde{A}).
\label{ferr1}
\end{equation}

Note that from \eqref{ferr1}, we have that for SPAI-GMRES-IR, the quantity $u_s\Vert E_i\Vert_\infty$ in \eqref{ass1} is on the order $f(n,k) \left(u_g + u_p \kappa_\infty(A)\right) \kappa_\infty(\tilde{A})$, and by \eqref{begmres}, the quantities $c_1$ and $c_2$ in \eqref{ass2} will be similar to those stated for five-precision GMRES-IR in \cite{h:21}. The difference will lie in the values of $\kappa_\infty(\tilde{A})$ and $\Vert \tilde{A}\Vert_F$, which here will depend on the $\varepsilon$ used in constructing the SPAI preconditioner. 

Using the estimate for the condition number of the preconditioned coefficient matrix in \eqref{kappaAtest}, we can thus give the rough bound
\begin{equation}
\frac{\Vert \dhat_i - d_i\Vert_\infty}{\Vert d_i\Vert_\infty} \lesssim f(n,k) \left(u_g + u_p \kappa_\infty(A)\right)  (1+2n\varepsilon)^2.
\label{ferr}
\end{equation}

Assume we are in the setting where $u=u_p=u_g$, which is commonly used in practice. 
Combining \eqref{ferr} with \eqref{convcondloose}, we thus must have roughly
\[
u_f\cond_2(A^T) \lesssim \varepsilon \lesssim u^{-1/2}\kappa_\infty(A)^{-1/2}
\]
in order to guarantee both that the SPAI algorithm will complete and that the GMRES-based iterative refinement scheme will converge.

We note that the bound \eqref{nrbe} on the normwise relative backward error that can be achieved with GMRES in precisions $u_g$ and $u_p$ assumes that we do not stop the iterations early. In practice, we will use some stopping criterion $\tau$ within GMRES. If $\tau$ is based on the backward error of the preconditioned system, then \eqref{nrbe} becomes
\begin{equation}
\frac{\Vert z_i - \tilde{A}\dhat_i\Vert_\infty}{\Vert\tilde{A}\Vert_\infty \Vert \dhat_i\Vert_\infty + \Vert z_i \Vert_\infty} \lesssim \max \left(f(n,k) \left(u_g + u_p \kappa_\infty(A)\right), \tau \right).
\label{nrbetau}
\end{equation}
In theory, we should then set $\tau$ to be on the same order as the other term in the maximum above if we do not want to affect the convergence conditions. In practice, however, the convergence conditions are often too strict, and thus we can use a larger $\tau$ and still maintain convergence of the iterative refinement process, perhaps with a larger $\phi_i$. In practice, choosing $\tau$ is largely done based on empirical observations, and is a tradeoff between the number of GMRES iterations per refinement step and the overall number of refinement steps; see, e.g., the discussion and experiments in \cite[Section 4.2]{h:21}.

Comparing the expected behavior of SPAI-GMRES-IR and GMRES-IR with full LU factorization, we expect that convergence for GMRES with the SPAI preconditioner will in general be slower than GMRES using the full LU factors, at least when we assume common choices of $\varepsilon$. However, the construction and application of the SPAI preconditioner may be less expensive and more parallelizable than for the full LU factors. This behavior and the resulting tradeoffs will be strongly problem dependent.

\section{Numerical Experiments}
\label{sec:experiments}

We now turn to an experimental evaluation of the numerical behavior of SPAI-GMRES-IR.
Our experiments are performed in MATLAB R2020a. In these experiments we use four precisions. For half precision, we use the \texttt{chop} library available at \texttt{https://github.com/higham/chop} and the associated functions available at \texttt{https://github.com/SrikaraPranesh/Multi\_precision\_NLA\_kernels}. For quadruple precision, we use the Advanpix Multiprecision Computing Toolbox \cite{adva-mct}. For single and double precision we use the built-in MATLAB datatypes. The code used to generate plots in this work can be found at \texttt{https://github.com/Noaman67khan/SPAI-GMRES-IR}.

The matrices tested come from the SuiteSparse Matrix Collection \cite{davis2011university}. We list tested matrices and their relevant properties in Table \ref{tab:mats}. 
In all tests, the right-hand sides are set to the vector with equal components and unit $2$-norm. For the convergence tolerance in GMRES, we use $\tau = 10^{-4}$ when the working precision is single and $\tau = 10^{-8}$ when the working precision is double, which responds to roughly the square root of the working precision. These are the default values which have been used in previous works, e.g., \cite{h:21}, \cite{carson2018accelerating}, although as mentioned in Section \ref{sec:spaigmresir}, this parameter must be tuned in practical applications. In Table \ref{tab:mats}, we split the matrices into two sets, those used with precisions $(u_f, u, u_r)=$ (single, double, quad), and those used with precisions $(u_f, u, u_r)=$  (half, single double), based on whether GMRES-IR with LU preconditioning converges with the given precisions and value of $\tau$.

In all cases, we set the initial sparsity pattern for SPAI to be the identity matrix. As noted by Sedlacek \cite[Section 3.1.2]{sedlacek2012sparse}, this can result in SPAI preconditioners with zero columns in the case that $A$ has a zero entry on the diagonal. We therefore restrict ourselves to problems with nonzero diagonal entries, but note that this could be remedied by either permuting $A$ or using the initial sparsity pattern of $A$, which, when SPAI is run on $A^T$, guarantees that we obtain a $\Pc$ with nonzero rows \cite[Theorem 3.1]{sedlacek2012sparse}.

For SPAI-GMRES-IR, we preprocess the matrices via a column scaling such that the largest element in absolute value in every column of $A^T$ is 1. A one-sided scaling was proposed in the context of using low-precision QR factorization to solve least squares problems in \cite[Section 3.2]{carson2020three}, which is precisely the task we perform in SPAI construction. The goal of this scaling is to avoid overflow in the low precision computations. To be precise, we run SPAI on $A^T D$ to obtain $M$, where $D$ is the diagonal matrix that performs the scaling, and then set $P = M^T D$. For all tests, we use $\beta=8$, which is in the range suggested by Sedlacek  \cite{sedlacek2012sparse}. Since we want to allow the approximate inverse to fill in as much as necessary so that \eqref{eq:trueres} is satisfied, we allow $\alpha$ to be as large as $\lceil n/\beta\rceil$. 

As a point of comparison, we also provide results for GMRES-IR with full LU factorization (which in this section we refer to as LU-GMRES-IR) as well as GMRES-IR with no preconditioner. 
To enable a fair comparison, for LU-GMRES-IR, we first test multiple reorderings of the matrix and then choose the one that provides the least number of nonzeros in $L$ and $U$. The reorderings tested are reverse Cuthill-McKee (rcm), approximate minimum degree (amd), nested dissection (nds), and column count (col). For LU-GMRES-IR, if $u_f = $ half produces LU factors containing \texttt{NAN} or \texttt{INF}, we implement the scaling advised in \cite{higham2019squeezing}.

We note that for SPAI-GMRES-IR, we use the natural ordering; as noted in \cite[p. 12]{benzi2000orderings}, there is no expected benefit to reordering for the case of nonfactorized sparse approximate inverses. 
We also note that in all variants of GMRES-IR tested here, including SPAI-GMRES-IR, we use the variant in which $u_g=u_p=u$, which is commonly used in practice.

\begin{table}[t]
\centering
\caption{Matrices used in experiments along with their properties. All matrices come from the SuiteSparse collection \cite{davis2011university}. }
\label{tab:mats}
\begin{tabular}{|c|c|c|c|c|c|}
\hline
Name & $n$ & $nnz$ & $nnz(A^{-1})$& $\kappa_\infty(A)$ & $\text{cond}_2(A^T)$ \\ \hline
\texttt{pores\_3} & 532 & 3474 & 213712 & $1.2e\pl 06$ & $1.7e\pl 05$  \\ \hline
\texttt{steam1} & 240 & 2248 & 57599 & $3.1e\pl 07$ & $2.8e\pl 03$ \\ \hline
\texttt{steam3} & 80 & 314 & 6315& $7.6e\pl 10$ & $5.6e\pl 03$ \\ \hline
\texttt{saylr1} & 238 & 1128 & 56644 & $1.6e\pl 09$ & $5.2e\pl 05$ \\ \hline \hline
\texttt{bfwa782} & 782 & 7514 & 458839 & $6.8e\pl 03$ & $1.3e\pl 03$ \\ \hline
\texttt{cage5} & 37 & 233 & 1369 & $2.9e\pl 01$ & $7.5e\pl 00$ \\ \hline
\texttt{gre\_115} & 115 & 421 & 13225 & $1.4e\pl 02$ & $3.7e\pl 01$ \\ \hline
\texttt{orsreg\_1} & 2205 & 14133 & 4862025 & $7.0e\pl 03$ & $5.9e\pl 03$ \\ \hline
\texttt{sherman4} & 1104 & 3786 & 298674 & $3.1e\pl 03$ & $1.2e\pl 03$ \\ \hline
\end{tabular}
\end{table}

For each linear system and given combination of precisions, we run SPAI-GMRES-IR with various values of $\varepsilon$, LU-GMRES-IR, and GMRES-IR with no preconditioning.  We note that it is incredibly difficult to give an accurate comparison of cost among the different methods; the relative cost of a GMRES iteration in GMRES-IR without preconditioning and SPAI-GMRES-IR and LU-GMRES-IR will be different, since the (preconditioned) coefficient matrix has different structure and density. SPAI-based preconditioners may also be more inexpensive to apply (and more parallelizable), since their application involves only a single SpMV rather than two triangular solves. Here we will give a qualitative comparison of approaches in terms of preconditioner size and number of GMRES iterations, but we note that high-performance experiments on large-scale matrices are necessary future work in order to give a quantitative comparison.

The tables for each setup give the size of the preconditioner in terms of number of nonzeros, the resulting condition number of the preconditioned system, as well as information about the number of GMRES-IR refinement steps and GMRES iterations per refinement step. For LU-GMRES-IR, the size of the preconditioner is given as $nnz(L+U)$, where $\Pi(QAQ^T)=LU$ where $Q$ is the best reordering among those tested and $\Pi$ is a permutation arising from partial pivoting. The best reorderings are noted in the tables. We stress that quantitative conclusions about the relative performance of LU-GMRES-IR versus SPAI-GMRES-IR cannot be drawn from the size of the preconditioner. 

In the final column in the tables, the first number gives the total number of GMRES iterations summed over all refinement steps, and the following tuple gives in position $i$ the number of GMRES iterations in refinement step $i$. Note that the number of refinement steps is given by the number of elements in the tuple.

The best case scenario is that SPAI preconditioning produces a preconditioner that is sparser than the full LU factorization and still effective in reducing the total number of GMRES iterations versus the unpreconditioned case. We again stress, however, that a full picture of expected relative performance cannot be given by the size of the preconditioner and the number of iterations. 
In Section \ref{sec:sdq}, we give examples using precisions $(u_f, u, u_r)=$ (single, double, quad) and in Section \ref{sec:hsd}, we use $(u_f, u, u_r)=$ (half, single, double). Finally, in Section \ref{sec:ssd}, we redo the experiments in Section \ref{sec:hsd} but now using $(u_f, u, u_r)=$ (single, single, double). These results demonstrate that there is often not a significant advantage, either in terms of the size of the preconditioner or the total number of GMRES iterations, to using $u_f$ higher than that dictated by the constraint \eqref{convcondloose}. That is, the resulting SPAI preconditioner will satisfy \eqref{eq:AM1} regardless of whether half or single precision is used, and thus we do not expect a significant difference in the behavior of the iterative refinement process.

\subsection{Experiments with $(u_f, u, u_r)=$ (single, double, quad)}
\label{sec:sdq}

Experiments using the settings $(u_f, u, u_r)=$ (single, double, quad) are shown in Table \ref{tab:sdq}. 
First, we observe that it is clear that SPAI-GMRES-IR converges in every case. 
The potential of SPAI-GMRES-IR is perhaps best exhibited for the matrix \texttt{steam1}. Here, even with reordering the computed LU factors contain 14,133 nonzeros, compared to the SPAI preconditioner, which contains 1,140 nonzeros when $\varepsilon=0.2$ and $1,303$ nonzeros when $\varepsilon=0.1$ (note that this is fewer than the number of nonzeros in $A$). In both cases, SPAI-GMRES-IR requires only 14 total GMRES iterations to converge to double precision accuracy; while this more than required by LU-GMRES-IR, it is still significantly fewer than the case where no preconditioner is used, which requires 543 total GMRES iterations.

For the matrices \texttt{pores\_3}, \texttt{steam3}, and \texttt{saylr1}, the results are as expected. The SPAI preconditioner, depending on the $\varepsilon$ parameter used, can have significantly fewer nonzeros than the LU factors, and the total number of GMRES iterations required for SPAI-GMRES-IR is somewhere in between that of LU-GMRES-IR and GMRES-IR with no preconditioning. We again stress that it is difficult to make quantitative comparisons; although the SPAI preconditioner may have close to the same number of nonzeros as the LU factors, it may be less expensive to apply and more parallelizable. We also note that the tradeoff in choosing $\varepsilon$ is clear: a larger $\varepsilon$ can lead to a sparser preconditioner, but a greater total number of GMRES iterations required for the iterative refinement process to converge.

%%%%% pores3 sqd %%%%%%%%

\begin{table}[]
\caption{Comparison of SPAI-GMRES-IR for different $\varepsilon$ values with LU-GMRES-IR and GMRES-IR with no preconditioner for test matrices using $(u_{f},u,u_{r})=$ (single, double, quad).}
\begin{center}
     {\footnotesize{

\begin{tabular}{|c|c|c|c|c|}
\hline
& Preconditioner & $\kappa_\infty(\tilde{A})$ & Precond. $nnz$ & \begin{tabular}[c]{@{}c@{}}GMRES-IR  steps/iterations\end{tabular} \\ \hline
\parbox[t]{2mm}{\multirow{4}{*}{\rotatebox[origin=c]{90}{\texttt{pores\_3}}}} & SPAI, $\varepsilon=0.5$ & $6.6e\pl 03$ & 3560 &  223(110, 113) \\ 
& SPAI, $\varepsilon=0.4$ & $3.8e\pl 03$ & 4871 &  174(86, 88) \\ 
& Full LU (amd) & $1.0e\pl 00$ & 9706 & 4(2, 2)\\ 
& None & $1.2e\pl 06$ & 0 & 1314(417, 456, 441) \\ \hline
\parbox[t]{2mm}{\multirow{4}{*}{\rotatebox[origin=c]{90}{\texttt{steam1}}}} & SPAI, $\varepsilon=0.2$ & $1.5e\pl 00$ & 1140 &  14(7, 7) \\
& SPAI, $\varepsilon=0.1$ & $1.5e\pl 00$ & 1303 &  14(7, 7) \\ 
& Full LU (amd) & $1.9e\pl 00$ & 14133 & 2(2)\\ 
& None & $3.1e\pl 07$ & 0 & 543(158, 193, 192) \\ \hline
\parbox[t]{2mm}{\multirow{4}{*}{\rotatebox[origin=c]{90}{\texttt{steam3}}}} & SPAI, $\varepsilon=0.5$ & $3.9e\pl 00$ & 244 &  31(9, 12, 10) \\ 
& SPAI, $\varepsilon=0.1$ & $1.9e\pl 00$ & 403 & 17(5, 6, 6) \\ 
& Full LU (amd) & $1.1e\pl 00$ & 483 & 2(2)\\ 
& None & $7.6e\pl 10$ & 0 & 221(61, 80, 80) \\ \hline
\parbox[t]{2mm}{\multirow{4}{*}{\rotatebox[origin=c]{90}{\texttt{saylr1}}}} & SPAI, $\varepsilon=0.4$ & $1.9e\pl 04$ &  1932 &  195(64, 66, 65) \\ 
& SPAI, $\varepsilon=0.3$ & $7.5e\pl 03$ & 3405 & 89(44, 45) \\ 
& Full LU (amd) & $1.0e\pl 00$ & 3607 & 5(2, 3) \\ 
& None & $1.6e\pl 09$ & 0 & 658(214, 229, 215)
\\ \hline
\end{tabular}
}}
\end{center}
\label{tab:sdq}
\end{table}
%\end{figure}

\subsection{Experiments with $(u_f, u, u_r) = $ (half, single, double)}
\label{sec:hsd}

We now turn to examples using precisions $(u_f, u, u_r) = $ (half, single, double), the results for which are shown in Table \ref{tab:hsd}. We again note that we have limited ourselves to examples where the use of half precision does not cause over/underflow in the SPAI construction despite the scaling used. 

From Table \ref{tab:hsd}, we can see that SPAI-GMRES-IR indeed converges in every case, although again we can see that there are clear tradeoffs related to the choice of $\varepsilon$. When we choose a larger $\varepsilon$, the preconditioner becomes sparser, but also less effective, requiring a greater total number of GMRES iterations to converge. For example, for the linear system with the matrix \texttt{sherman4}, $\varepsilon=0.5$ leads to a very sparse preconditioner (nearly diagonal), but requires almost as many GMRES iterations across all refinement steps as the unpreconditioned case. 

For the problem with matrix \texttt{gre\_115}, we can see that for the choice of the smaller $\varepsilon$, $\varepsilon = 0.3$, the resulting SPAI preconditioner has a greater number of nonzeros than the LU preconditioner, and is also less effective in terms of reducing the total number of GMRES iterations. 

The matrix \texttt{orsreg\_1} perhaps best exhibits the potential of SPAI-GMRES-IR. Here, with $\varepsilon=0.3$, the SPAI preconditioner has over $12\times$ fewer nonzeros than the LU factors  and still reasonably improves the convergence rate ($5.8\times$ fewer GMRES iterations than the unpreconditioned case).

%%%%% pores3 sqd %%%%%%%%

\begin{table}[]
\caption{Comparison of SPAI-GMRES-IR for different $\varepsilon$ values with LU-GMRES-IR and GMRES-IR with no preconditioner for test matrices using $(u_{f},u,u_{r})=$ (half, single, double).}
\begin{center}
     {\footnotesize{

\begin{tabular}{|c|c|c|c|c|}
\hline
& Preconditioner & $\kappa_\infty(\tilde{A})$ & Precond. $nnz$ & \begin{tabular}[c]{@{}c@{}}GMRES-IR  steps/iterations\end{tabular} \\ \hline
%
% \parbox[t]{2mm}{\multirow{4}{*}{\rotatebox[origin=c]{90}{\texttt{pores\_3}}}} & SPAI, $\varepsilon=0.5$ & $7.3e\pl 03$ & 3235 &  222(111, 111) \\ 
% & SPAI, $\varepsilon=0.4$ & $7.1e\pl 03$ & 3594 &  195(97, 98) \\ 
% & Full LU (amd) & $8.1e\pl 12$ & 283024 & $-$\\ 
% & None & $1.2e\pl 06$ & 0 & $-$ \\ \hline
% %
% \parbox[t]{2mm}{\multirow{4}{*}{\rotatebox[origin=c]{90}{\texttt{steam1}}}} & SPAI, $\varepsilon=0.2$ & Inf & 220 &  $-$ \\
% & SPAI, $\varepsilon=0.1$ & Inf & 220 &  $-$ \\ 
% & Full LU (amd) & $1.3e\pl 14$ & 57600 & $-$\\ 
% & None & $3.1e\pl 07$ & 0 & 2450(50, 240,  240,\ldots,  240) \\ \hline
% % %
% % \parbox[t]{2mm}{\multirow{4}{*}{\rotatebox[origin=c]{90}{\texttt{steam3}}}} & SPAI, $\varepsilon=0.5$ & $3.9e\pl 00$ & 244 &  31(9, 12, 10) \\ 
% % & SPAI, $\varepsilon=0.1$ & $1.9e\pl 00$ & 403 & 17(5, 6, 6) \\ 
% % & Full LU (amd) & $1.1e\pl 00$ & 483 & 2(2)\\ 
% % & None & $7.6e\pl 10$ & 0 & 221(61, 80, 80) \\ \hline
% % %
% % \parbox[t]{2mm}{\multirow{4}{*}{\rotatebox[origin=c]{90}{\texttt{saylr1}}}} & SPAI, $\varepsilon=0.4$ & $1.9e\pl 04$ &  1932 &  195(64, 66, 65) \\ 
% % & SPAI, $\varepsilon=0.3$ & $7.5e\pl 03$ & 3405 & 89(44, 45) \\ 
% % & Full LU (amd) & $1.0e\pl 00$ & 3607 & 5(2, 3) \\ 
% % & None & $1.6e\pl 09$ & 0 & 658(214, 229, 215)
% % \\ \hline
% %
%
\parbox[t]{2mm}{\multirow{4}{*}{\rotatebox[origin=c]{90}{\texttt{bfwa782}}}} &  SPAI, $\varepsilon=0.5$ & $1.1e\pl 03$ & 6271 & 164(75, 89) \\ 
  & SPAI, $\varepsilon=0.3$ & $5.0e\pl 02$ & 11430 & 114(54, 60) \\ 
  & Full LU (amd) & $2.1e\pl 00$ & 21838 & 7(3, 4) \\ 
  & None & $6.8e\pl 03$ & 0 & 381(172, 209)\\ \hline
\parbox[t]{2mm}{\multirow{4}{*}{\rotatebox[origin=c]{90}{\texttt{cage5}}}} & SPAI, $\varepsilon=0.5$ & $9.9e\pl 00$ &  101 &  16(8, 8) \\ 
& SPAI, $\varepsilon=0.3$ & $3.9e\pl 00$ &  213 & 12(6, 6) \\ 
& Full LU (amd) & $1.0e\pl 00$ & 359 & 2(2) \\ 
& None & $2.9e\pl 01$ & 0 & 25(13, 12) \\ \hline
\parbox[t]{2mm}{\multirow{4}{*}{\rotatebox[origin=c]{90}{\texttt{gre\_115}}}} &   SPAI, $\varepsilon=0.5$ & $5.8e\pl 02$ &  725 &  48(24, 24) \\ 
  & SPAI, $\varepsilon=0.3$ & $1.8e\pl 01$ & 1719 & 21(10, 11) \\ 
  & Full LU (nds) & $1.0e\pl 00$ & 1551 & 2(2) \\
  & None & $1.4e\pl 02$ & 0 & 100(49, 51) \\ \hline
 \parbox[t]{2mm}{\multirow{4}{*}{\rotatebox[origin=c]{90}{\texttt{orsreg\_1}}}} &
   SPAI, $\varepsilon=0.5$ & $1.7e\pl 02$ &  9261 &  108(29, 45, 34) \\ 
  & SPAI, $\varepsilon=0.3$ & $1.3e\pl 02$ & 11120 & 61(23, 38) \\
  & Full LU (rcm) & $2.2e\pl 00$ & 133634 & 9(4, 5) \\ 
  & None & $7.0e\pl 03$ & 0 & 352(107, 150, 95) \\ \hline
\parbox[t]{2mm}{\multirow{4}{*}{\rotatebox[origin=c]{90}{\texttt{sherman4}}}} &   SPAI, $\varepsilon=0.5$ & $1.6e\pl 03$ &  1386 & 140(67, 73) \\ 
 & SPAI, $\varepsilon=0.3$ & $5.0e\pl 02$ & 8496 & 74(35, 39) \\ 
 & Full LU (amd) & $1.8e\pl 00$ & 14211 & 5(2, 3) \\ 
 & None & $3.1e\pl 03$ & 0 & 178(85, 93) \\ \hline
\end{tabular}
}}
\end{center}
\label{tab:hsd}
\end{table}
%\end{figure}

\subsection{ Experiments with $(u_f, u, u_r) = $ (single, single, double)}
\label{sec:ssd}

We now test the same problems used in Section \ref{sec:hsd}, now using single precision instead of half precision for $u_f$. The results are presented in Table \ref{tab:ssd}. The analysis in Section \ref{sec:analysis} indicates that as long as, say, half precision is sufficient for producing a preconditioner satisfying \eqref{convcondloose}, there is in general no expected benefit to using higher precision in the SPAI construction in terms of preconditioner quality. We do recall, however, that the number of nonzeros in the preconditioner may be different when higher versus lower precision is used. For easy comparison, in Table \ref{tab:ssd}, we give in red bracketed format the corresponding condition number, preconditioner size, and total number of GMRES iterations from Table \ref{tab:hsd}.

For all test problems, notice that the results for SPAI-GMRES-IR  either do not change at all or change only very slightly. In cases where there is a slight difference, note that it is not necessarily clear that $u_f$ = single or $u_f=$ half is better. For example, for the problem with \texttt{gre\_115} and $\varepsilon=0.5$, $u_f=$ half leads to 48 total GMRES iterations whereas $u_f=$ single leads to 51. Conversely, for the problem with \texttt{orsreg\_1} and $\varepsilon=0.5$, $u_f=$ half leads to 108 total GMRES iterations whereas $u_f =$ single leads to 97. We note that for both these cases, the SPAI preconditioner has the same number of nonzeros regardless of whether $u_f=$ half or $u_f=$ single is used.

It is also not necessarily clear that $u_f$ = single or $u_f=$ half is better in terms of the size of the preconditioner. For example, for the problem with \texttt{bfwa782} and $\varepsilon=0.3$, the SPAI preconditioner has 11,430 nonzeros when $u_f=$ half, but 11,470 when $u_f=$ single. For the problem with \texttt{orsreg\_1} and $\varepsilon=0.3$, the SPAI preconditioner has 11,120 nonzeros when $u_f=$ half, and 11,025 when $u_f=$ single. 

In any case, we can see that this data supports our theoretical conclusions: in general, there is no expected benefit to using a higher precision than that dictated by \eqref{convcondloose} for SPAI in terms of preconditioner quality. Note that this is in sharp contrast to the case of LU preconditioning, where the precision used for the LU factorization makes a significant difference in preconditioner quality. This intuitively makes sense since LU factorization is a direct algorithm rather than an iterative one like SPAI that iterates until the stopping criterion is met. 

One interesting anomaly that we point out is that from Table \ref{tab:ssd}, we can see that the LU factors can have significantly fewer nonzeros when we use $u_f=$ half versus $u_f=$ single. This difference is most evident for the problem with matrix \texttt{orsreg\_1}, where $u_f=$ half results in preconditioners with 133,634 nonzeros and $u_f=$ single results in preconditioners with 330,910 nonzeros. We expect this is due to underflow; the true LU factors likely contain many small elements which become zeros when half precision is used.

%%%%% pores3 sqd %%%%%%%%

\begin{table}[]
\caption{Comparison of SPAI-GMRES-IR for different $\varepsilon$ values with LU-GMRES-IR and GMRES-IR with no preconditioner for test matrices using $(u_{f},u,u_{r})=$ (single, single, double).}
\begin{center}
    {\footnotesize{

\begin{tabular}{|c|c|c|c|c|}
\hline
& Preconditioner & $\kappa_\infty(\tilde{A})$ & Precond. $nnz$ & \begin{tabular}[c]{@{}c@{}}GMRES-IR  steps/iterations\end{tabular} \\ \hline
\parbox[t]{2mm}{\multirow{4}{*}{\rotatebox[origin=c]{90}{\texttt{bfwa782}}}} &    SPAI, $\varepsilon=0.5$ & $1.1e\pl 03$ \cred{[$1.1e\pl 03$]} & 6261 \cred{[6271]} & 166(74, 92) \cred{[164]} \\ 
  & SPAI, $\varepsilon=0.3$ & $5.0e\pl 02$ \cred{[$5.0e\pl 02$]} & 11470 \cred{[11430]} & 114(54, 60) \cred{[114]} \\
  & Full LU (amd) & $1.0e\pl 00$ \cred{[$2.1e\pl 00$]} & 21848 \cred{[21838]} & 1(1) \cred{[7]}\\ 
    & None & $6.8e\pl 03$ & 0 & 381(172, 209)\\ \hline
\parbox[t]{2mm}{\multirow{4}{*}{\rotatebox[origin=c]{90}{\texttt{cage5}}}} & SPAI, $\varepsilon=0.5$ & $9.9e\pl 00$ \cred{[$9.9e\pl 00$]} &  101 \cred{[101]} & 16(8, 8) \cred{[16]}\\ 
& SPAI, $\varepsilon=0.3$ & $3.9e\pl 00$ \cred{[$3.9e\pl 00$]} &  213 \cred{[213]}& 12(6, 6) \cred{[12]}\\ 
& Full LU (amd) & $1.0e\pl 00$ \cred{[$1.0e\pl 00$]} & 359 \cred{[359]}& 1(1) \cred{[2]}\\ 
& None & $2.9e\pl 01$ & 0 & 25(13, 12) \\ \hline
\parbox[t]{2mm}{\multirow{4}{*}{\rotatebox[origin=c]{90}{\texttt{gre\_115}}}} &   SPAI, $\varepsilon=0.5$ & $6.2e\pl 02$ \cred{[$5.8e\pl 02$]} &  725 \cred{[725]}&  51(24, 27) \cred{[48]}\\ 
 &  SPAI, $\varepsilon=0.3$ & $1.7e\pl 01$ \cred{[$1.8e\pl 01$]} & 1739 \cred{[1719]} & 20(10, 10) \cred{[21]}\\ 
 & Full LU (nds) & $1.0e\pl 00$ \cred{[$1.0e\pl 00$]} & 1556 \cred{[1551]}& 1(1) \cred{[2]} \\ 
  & None & $1.4e\pl 02$ & 0 & 100(49, 51) \\ \hline
 \parbox[t]{2mm}{\multirow{4}{*}{\rotatebox[origin=c]{90}{\texttt{orsreg\_1}}}} &
  SPAI, $\varepsilon=0.5$ & $1.4e\pl 02$ \cred{[$1.7e\pl 02$]}&  9261 \cred{[9261]} &  97(25, 40, 32) \cred{[108]} \\ 
  & SPAI, $\varepsilon=0.3$ & $1.1e\pl 02$ \cred{[$1.3e\pl 02$]}& 11025 \cred{[11120]} & 60(22, 38) \cred{[61]} \\ 
  & Full LU (rcm) & $1.0e\pl 00$ \cred{[$2.2e\pl 00$]} & 330910 \cred{[133634]} & 1(1) \cred{[9]} \\ 
   & None & $7.0e\pl 03$ & 0 & 352(107, 150, 95) \\ \hline
\parbox[t]{2mm}{\multirow{4}{*}{\rotatebox[origin=c]{90}{\texttt{sherman4}}}} & SPAI, $\varepsilon=0.5$ & $1.6e\pl 03$ \cred{[$1.6e\pl 03$]} & 1385 \cred{[1386]}& 140(67, 73) \cred{[140]} \\ 
& SPAI, $\varepsilon=0.3$ & $5.0e\pl 02$ \cred{[$5.0e\pl 02$]} & 8499 \cred{[8496]}& 74(35, 39)\cred{[74]} \\
& Full LU (amd) & $1.0e\pl 00$ \cred{[$1.8e\pl 00$]} & 14211 \cred{[14211]} & 1(1) \cred{[5]}\\ 
 & None & $3.1e\pl 03$ & 0 & 178(85, 93) \\ \hline
\end{tabular}}}
\end{center}
\label{tab:ssd}
\end{table}
%\end{figure}

\section{Conclusions and future work}
\label{sec:conclusion}

In this work we explored the use of sparse approximate inverse preconditioners within mixed precision GMRES-based iterative refinement as a method of solving sparse linear systems. 

We first analyzed the computation of an SPAI preconditioner in some finite precision $u_f$. The stopping criterion used in constructing each column of the sparse approximate inverse $M$ is based on the size of the measured residual norm $e_k-A^Tm_k$ being less than some user-specified parameter $\varepsilon$. Thus for a given $\varepsilon$, our analysis shows that we must use $u_f \lesssim \varepsilon \cond_2^{-1}(A^T)$ in order to guarantee that this stopping criterion can be achieved. As long as this constraint on $u_f$ holds, then the computed $\Mhat$ satisfies the bounds proved for the exact arithmetic case in \cite{grote1997parallel} up to a constant factor. An interesting point is that there is no significant benefit in terms of resulting preconditioner quality to using a precision higher than dictated by this level.

We then prove conditions under which GMRES left preconditioned by the SPAI preconditioner will be backward stable with respect to the preconditioned system. As long as the constraint $u_f \lesssim \varepsilon \cond_2^{-1}(A^T)$ is satisfied, meaning that an approximate preconditioner with parameter $\varepsilon$ is computable in precision $u_f$, then the constraints for convergence of five-precision SPAI-GMRES-IR are essentially the same as those for the five-precision GMRES-IR given in \cite{h:21}. Our numerical experiments confirm the expected behavior of SPAI-GMRES-IR.

While we have started with one variant of sparse approximate inverses here, it is clear that other approximate preconditioners, such as incomplete LU factorization, factorized sparse approximate inverses, or sparse approximate inverses based on a fixed sparsity pattern may be suitable to use within mixed precision Krylov-subspace-based iterative refinement schemes. In particular, we believe that extension of the analysis here to modified SPAI (see, e.g., \cite{huckle2007frobenius}) as well as factorized SPAI (see, e.g.,  \cite{kolotilina1993factorized}, \cite{yeremin2004factorized}) should be straightforward.

\bibliographystyle{siamplain}
\bibliography{spai}
\end{document}